\documentclass[11pt]{amsart} 

\usepackage[all]{xy}
\usepackage{amsthm}
\usepackage{mathrsfs}
\usepackage{amssymb}
\usepackage{amsmath}
\usepackage{array, url}
\usepackage{stmaryrd}
\usepackage{mathtools}
\usepackage{wasysym}
\usepackage{tensor}
\usepackage{mathdots}
\usepackage{amsbsy}
\usepackage{tikz-cd}

\usepackage[T2A,T1]{fontenc}
\DeclareSymbolFont{cyrillic}{T2A}{cmr}{m}{n}
\DeclareMathSymbol{\Sha}{\mathalpha}{cyrillic}{216}

\usepackage{latexsym}
\usepackage{amsxtra}

\let\oldmarginpar\marginpar
\renewcommand\marginpar[1]{\-\oldmarginpar[\raggedleft\footnotesize #1]%
{\raggedright\footnotesize #1}}

\newcounter{firstnumber}[section]
\newcounter{secondnumber}[firstnumber]
\newcounter{thirdnumber}[secondnumber]
\newcounter{fourthnumber}[thirdnumber]
\newcounter{fifthnumber}[fourthnumber]
\newcounter{currentdepth}

\renewcommand{\thefirstnumber}{\arabic{section}.\arabic{firstnumber}}
\renewcommand{\thesecondnumber}{\thefirstnumber.\arabic{secondnumber}}

\newcommand{\segment}[2]{%
    \setcounter{currentdepth}{1}%
    \def\thesubsection{\thefirstnumber}%
    \refstepcounter{firstnumber}\label{#1}%
    \addtocounter{subsection}{-1}%
    \subsection{#2}}
\newcommand{\ssegment}[2]{%
    \setcounter{currentdepth}{2}%
    \def\thesubsection{\thesecondnumber}%
    \refstepcounter{secondnumber}\label{#1}%
    \addtocounter{subsection}{-1}%
    \subsection{#2}}

\newcommand{\bs}{\backslash}

\newcommand{\set}[2]{\big\{ #1 \; \big| \; #2 \big\} }

\newcommand{\xto}{\xrightarrow}

\newcommand{\surj}{\twoheadrightarrow}

\renewcommand{\Im}{\operatorname{Im}}
\newcommand{\ord}{\operatorname{ord}}
\newcommand{\rec}{\operatorname{rec}}
\newcommand{\rank}{\operatorname{rank}}

\newcommand{\Spec}{\operatorname{Spec}}

\newcommand{\Pic}{\operatorname{Pic}}

\newcommand{\Li}{\operatorname{Li}}

\newcommand{\Aut}{\mathrm{Aut}\,}

\newcommand{\m}[1]{\mathrm{#1}}

\newcommand{\bb}[1]{\mathbb{#1}}

\newcommand{\la}{\lambda}

\newcommand{\ze}{\zeta}
\newcommand{\ga}{\gamma}
\newcommand{\al}{\alpha}
\newcommand{\be}{\beta}

\newcommand{\Om}{\Omega}

\newcommand{\Gm}{{\mathbb{G}_m}}
\newcommand{\Qp}{{\QQ_p}}
\newcommand{\Zp}{{\ZZ_p}}

\newcommand{\Qbar}{{\overline{\QQ}}}

\newcommand{\ZZ}{\bb{Z}}
\newcommand{\CC}{\bb{C}}
\newcommand{\GG}{\mathbb{G}}

\newcommand{\QQ}{\bb{Q}}

\newcommand{\PP}{\bb{P}}

\newcommand{\Dd}{\mathcal{D}}
\newcommand{\Tt}{\mathcal{T}}

\renewcommand{\AA}{\bb{A}}

\newcommand{\Oo}{\mathcal{O}}
\newcommand{\Gg}{\mathcal{G}}

\newcommand{\Hh}{\mathcal{H}}

\newcommand{\Xx}{\mathcal{X}}
\newcommand{\Yy}{\mathcal{Y}}

\newcommand{\Ee}{\mathcal{E}}
\newcommand{\Ll}{\mathcal{L}}

\newcommand{\inv}{^{-1}}

\newcommand{\areq}{\ar@{=}}
\newcommand{\suphook}{\ar@{^(->}}
\newcommand{\subhook}{\ar@{_(->}}

\newcommand{\inj}{\hookrightarrow}

\newcommand{\loc}{\operatorname{loc}}

\newcommand{\g}{\gamma}
\newcommand{\G}{\Gamma}
\newcommand{\rTo}{\to}

\newcommand{\Sel}{\operatorname{Sel}}

\newcommand{\et}{{\textrm {\'et}}}

\newcommand{\ZSinv}{\ZZ[S\inv]}

\newcommand{\bX}{\bar X}
\newcommand{\bF}{\bar F}
\newcommand{\pip}{\pi_{[2]}}
\newcommand{\Z}{\bb Z}
\newcommand{\cH}{\Hh}
\newcommand{\Q}{\QQ}
\newcommand{\cG}{\Gg}
\renewcommand{\t}{\tau}

\newcommand{\tbX}{\tilde \bX}
\newcommand{\ra}{\to}
\newcommand{\rOnto}{\surj}

\renewcommand{\O}{\Oo}
\newcommand{\cX}{\Xx}
\newcommand{\A}{\AA}
\newcommand{\z}{\zeta}
\renewcommand{\P}{\PP}
\newcommand{\cE}{\Ee}
\newcommand{\D}{\Delta}
\renewcommand{\a}{\al}
\renewcommand{\b}{\be}
\newcommand{\Exp}{\operatorname{Exp}}
\newcommand{\bE}{\bar E}
\newcommand{\piet}{\pi_1^\et(\bar X, b)}

\newcommand{\Fr}{\operatorname{Fr}}

\newcommand{\Id}{\operatorname{Id}}

\begin{document}

\SelectTips{cm}{11}

\title[A non-abelian conjecture]{A non-abelian conjecture of Tate--Shafarevich type for hyperbolic curves}

\author[Balakrishnan]{Jennifer S. Balakrishnan}
\author[Dan-Cohen]{Ishai Dan-Cohen}
\author[Kim]{Minhyong Kim}
\author[Wewers]{Stefan Wewers}

\date{\today}

\begin{abstract}
Let $X$ denote a hyperbolic curve over $\QQ$ and let $p$ denote a prime of good reduction. The third author's approach to integral points, introduced in \cite{kimi} and \cite{kimii}, endows $X(\Zp)$ with a nested sequence of subsets $X(\Zp)_n$ which contain $X(\ZZ)$. These sets have been computed in a range of special cases \cite{KimMassey, BalakAppendix, CKtwo, mtmue}; there is good reason to believe them to be practically computable in general. In 2012, the third author announced the conjecture that for $n$ sufficiently large, $X(\ZZ) = X(\Zp)_n$. This conjecture may be seen as a sort of compromise between the abelian confines of the BSD conjecture and the profinite world of the Grothendieck section conjecture. After stating the conjecture and explaining its relationship to these other conjectures, we explore a range of special cases in which the new conjecture can be verified.


\end{abstract}


\maketitle

\section{Introduction}

\segment{d41}{}
When $E/\QQ$ is an elliptic curve, the conjecture of Birch and Swinnerton-Dyer predicted the following phenomenon:
\[
 L(E,1)\neq 0 \Rightarrow | E(\QQ) | < \infty.
\]
This is now a theorem, strikingly realized 
by the  process of annihilating the Mordell-Weil group
with the $L$-value in question \cite{Kolyvagin, Kato}.
When we move to the realm of   hyperbolic curves, that is, curves with non-abelian
geometric fundamental groups, we have suggested elsewhere an extension of this
 connection between Diophantine finiteness
and non-vanishing of $L$-values \cite{CoatesKim, KimpL}, even though it has thus far proved
difficult to formulate it in precise terms.

\segment{d42}{}
The goal of this paper is to extend a different part of the constellation of conjectures surrounding BSD, namely, the finiteness of the Tate-Shafarevich group $\Sha_E$. To explain this, we begin by turning our attention to a different conjecture, namely Grothendieck's section conjecture. Let $X$ be a compact hyperbolic curve over $\QQ$, let $\bar X$ denote the base change of $X$ to an algebraic closure of $\QQ$ with Galois group $G$, and let $b$ be a $\QQ$-valued point of $X$. Then according to the conjecture, the map
\[
x\mapsto [\pi_1^\et(\bX;b,x)]
\]
that associates to a rational point $x$ the $\piet$-torsor of paths from $b$ to $x$ defines a bijection
\[
X(\Q)\simeq H^1(G, \piet).
\]
Returning to the special case of an elliptic curve $E$, our point of departure is the apparent similarity between this bijection and a certain isomorphism implied by the conjectured finiteness of $\Sha_E$, namely
\[
\tag{*}
E(\Q)\otimes_{\Z} \Q_p\simeq H^1_{\Z}(G, H_1(\bE, \Q_p)).
\]
Here, the subscript `$\Z$' refers to the cohomology classes for the group
$G$ that are {\em crystalline} at $p$, and zero at all $v\neq p$.

\segment{d43}{}
For the conjecture being presented here, we let $\Xx \to \Spec \ZZ$ be a
 \emph{regular minimal $\ZZ$-model of a hyperbolic curve} (see \ref{13a} for a precise definition); its generic fiber $X = \Xx_\QQ$ need not be proper. We let $b$ be an integral base point (possibly tangential), and assume $p$ is a prime of good reduction for $\Xx$ and $b$. Between the profinite fundamental group of the section conjecture, and the first \'etale homology of segment \ref{d41}, equation (*), lies the \emph{unipotent $p$-adic \'etale fundamental group} $U$ of $X_{\bar \QQ}$ at $b$. We let $U_n$ denote its $n^\m{th}$ quotient along the descending central series. Let $G_\QQ$ denote the total Galois group of $\QQ$ and let $G_p$ denote the total Galois group of $\Qp$. Following \cite{kimi, kimii}, we consider the subspace
 \[
 H^1_f(G_p, U_n) \subset H^1(G_p, U_n)
 \] 
consisting of $G_p$-equivariant $U_n$-torsors which are \emph{crystalline}. We also consider a certain subspace
\[
\Sel^n(\Xx) \subset H^1(G_{\QQ}, U_n), 
\]
the \emph{Selmer scheme of $\Xx$}; roughly speaking, it parametrizes those torsors which are crystalline at $p$ and in the image of $\Xx(\ZZ_v)$ (we say \textit{locally geometric}) for $v\neq p$. For each $n$ these fit into a commuting square like so,
\[
\xymatrix{
X(\ZZ) \ar[d]_-j \ar[r] &
X(\Zp) \ar[d]^-{j_p} \\
\Sel^n(\Xx) \ar[r]_-{\loc_p} &
H^1_f(G_p, U_n) }
\]
and we define
\[
\Xx(\Zp)_n := j_p\inv\big( \loc_p (\Sel^n(\Xx)) \big).
\]
These form a nested sequence of subsets like so.
\[
\Xx(\Zp) \supset \Xx(\Zp)_1 \supset \Xx(\Zp)_2 \supset \cdots \supset \Xx(\ZZ).
\]
The conjecture, which was first proposed by M.K. in his lectures at the I.H.E.S. in February of 2012, is as follows.

\medskip \noindent
\textbf{Conjecture} (\ref{13i} below)\textbf{.}
Equality $\Xx(\Zp)_n = \Xx(\ZZ)$ is obtained for large $n$.

\segment{}{}
We also suggest a variant of the Selmer scheme $\Sel_S^n(\Xx)$ of $\Xx$, suited to computing the $\ZSinv$-valued points of $\Xx$ for $S$ a finite set of primes, by dropping the local geometricity condition over $S$. This gives rise to a square
\[
\xymatrix{
X(\ZSinv) \ar[d]_-{j^S} \ar[r] &
X(\Zp) \ar[d]^-{j_p} \\
\Sel^n_S(\Xx) \ar[r]_-{\loc_p} &
H^1_f(G_p, U_n), }
\]
and to an associated sequence of subsets
\[
\Xx(\Zp) \supset \Xx(\Zp)_{S,1} \supset \Xx(\Zp)_{S,2} \supset \cdots \supset \Xx(\ZZ),
\]
for which equality $\Xx(\Zp)_{S,n} = \Xx(\ZSinv)$ may hold for large $n$.

\segment{}{}
These constructions are based on the third author's approach to integral points, introduced in \cite{kimi} and \cite{kimii}. The relationship to the section conjecture has been explored before. For instance in \cite{KimRemarks}, the third author shows that if $\Xx(\Zp)_n \neq \Xx(\Zp)$ for some $n$, then the section conjecture would in principle allow one to obtain a computable bound on the height of rational points. Our present conjecture, however, is quite different in flavor from the section conjecture and its direct consequences. It shares more with the conjectures of Tate--Shafarevich and Birch--Swinnerton-Dyer, both in terms of concreteness and in terms of computability. Indeed, like the BSD conjecture, the present conjecture can actually be tested numerically.

\segment{}{}
Our principal goal below is to do just that. Work completed elsewhere allows us to verify our conjecture in a range of cases. New in this article is the case of a punctured elliptic curve of rank zero: we are able to compute the sets $\Xx(\Zp)_2$, and subsequently to verify the conjecture for many such curves. This computation is based on a study of the \emph{unipotent Kummer map}
\[
j_v: X(F_v) \to H^1(G_v, U_2)
\]
for a punctured elliptic curve $X$ over a local field $F_v$ of residue characteristic $v \neq p$. Our main theorem (\ref{Ner0}) says that the $p$-adic height function can be retrieved from this map. This is of interest in its own right, and is suggestive of the possibility of obtaining functions on the local points from higher quotients of the unipotent fundamental group which might play a role similar to the role played by heights here. This point of view is implicit in the third author's work on nonabelian reciprocity laws \cite{KimRecip} and in ongoing joint work between him and Jonathan Pridham.

\segment{}{}
As explained in \cite{kimi, kimii}, the key to computing the map $j_p$ is its equivalence with a certain $p$-adic analog 
\[
\al: \Xx(\Zp) \to U_n^{DR}/F^0
\]
of the \emph{higher Albanese map} of Richard Hain \cite{HainHigher} through a lifting of the Bloch-Kato exponential to the unipotent level obtained via the unipotent $p$-adic Hodge theory of Martin Olsson \cite{OlssonTowards}. The $p$-adic unipotent Albanese map $\al$ is given in coordinates by certain $p$-adic iterated integrals (known also as Coleman functions), and there are fairly well-established methods for producing explicit formulas for the resulting iterated integrals on the one hand, and for computing $p$-adic approximations of their values on the other. For instance, the case of the thrice punctured line was treated by Furusho \cite{FurushoI, FurushoII} and by Besser--de Jeu \cite{Lip}. The problem of explicit determination of the unipotent Albanese map for punctured elliptic curves in depth two is treated by Kim \cite{KimMassey} and Balakrishnan--Kedlaya--Kim \cite{BalakAppendix}. The problem of computing Coleman functions on hyperelliptic curves is treated by Balakrishnan--Bradshaw--Kedlaya \cite{BalakBK} and by Balakrishnan \cite{BalakIterated}. 

\segment{}{}
We turn to the map
\[
\loc_p:\Sel^n(\Xx) \to H^1_f(G_p, U_n).
\]
This is actually an algebraic map of finite-type affine $\Qp$-schemes; its target is in fact isomorphic to affine space. Let $\Ll(n)$ denote the ideal defining its scheme-theoretic image. As explained in \cite{kimii}, as soon as $\Ll(n) \neq 0$, $\Xx(\Zp)_n$  becomes finite. Moreover, several well known motivic conjectures (Fontaine-Mazur-Jannsen, Bloch-Kato) imply that for $n$ large, 
\[
j_p^*\Ll(n+1) \supsetneq j_p^*\Ll(n), 
\]
and, in fact, the larger ideal contains elements that are algebraically independent of the elements in $j^*\Ll(n)$.\footnote{Technically speaking, the pullback $j_p^*$ appearing here may be thought of as a pullback of locally analytic functions on associated $p$-adic analytic spaces.} So a point in the common zero set for all $n$ should be there for a good reason; our conjecture expresses the belief that such a point must belong to $\Xx(\ZZ)$.

\segment{}{}
The study of the ideals $\Ll(n)$ relates not only the the plausibility of our conjecture, but also to its usefulness. Explicit computation of these ideals has been achieved in a range of special cases. An approach to the case of the thrice punctured line using the methods of mixed Tate motives is currently under development in a sequence of articles by Dan-Cohen and Wewers \cite{CKtwo, mtmue, mtmueII}. The case of punctured elliptic curves in depth two was treated by Kim \cite{KimMassey} and Balakrishnan--Kedlaya--Kim \cite{BalakAppendix}. The case of punctured hyperelliptic curves of genus equal to the Mordell-Weil rank of their Jacobian is treated in Balakrishnan--Besser--M\"uller \cite{BalBesMul}. The case of punctured elliptic curves of rank zero is treated in section \ref{S6} below. We believe strongly in the feasibility of computing the ideals $\Ll(n)$ and subsequently the loci $\Xx(\Zp)_n$ (as well as their $S$-integral variants) in a range of cases far beyond those mentioned above and detailed below. Such computations will provide powerful tools for bounding the number of ($S$-)integral points. If the conjecture holds, then bounds obtained in this way can be made sharp.

\segment{}{}
We begin in section \ref{S2} by giving a careful construction of $\Sel^n(\Xx)$. Our construction, which is a bit more elaborate than indicated above, relies on the work done in \cite{kimii} to endow $\Sel^n(\Xx)$ with the structure of an affine, finite-type $\Qp$-scheme for which the map $\loc_p$ is algebraic. In section \ref{S3}, after restating the conjecture, we discuss again in more detail its relationship to the finiteness of $\Sha$ and to the section conjecture, as well as the computability of the local Kummer map $j_p$ via the $p$-adic unipotent Albanese map.

\segment{}{}
The remainder of the article is devoted to discussing several special cases in which we are able to compute the loci $\Xx(\Zp)_{n,S}$ and so to obtain numerical evidence for the conjecture. Section \ref{S4} is devoted to proving a preliminary theorem to be used in our study of punctured elliptic curves of rank zero in section \ref{S6} below. 

Fix a prime $p \neq 2$. Let $F_v$ be a finite unramified extension of $\QQ_l$ for $l \neq p$ and let $E_v$ be an elliptic curve over $F_v$. We let $G_v$ denote the total Galois group of $F_v$. We fix a certain tangent vector $b$ at $O$ which serves as base point for the level $2$ quotient of the $p$-adic \'etale unipotent fundamental group $U_2$ of $X = E \setminus \{O\}$. Let
\[
\log \chi: G_v^\m{ab} \to \Qp
\]
denote the $p$-adic logarithm of the cyclotomic character. Let $j_v$ denote the local unipotent Kummer map
\[
X(F_v) \to H^1(G_v, U_2). 
\]
As we explain in segment \ref{d21}, the map
\[
H^1(G_v, \Qp(1)) \to H^1(G_v, U_2)
\]
induced by the inclusion $\Qp(1) \subset U_2$ is bijective. This allows us to regard $j_v$ as a map to $H^1(G_v, \Qp(1))$. Using the cup product 
\[
H^1(G_v, \Qp) \times H^1(G_v, \Qp(1)) \to H^2(G_v, \Qp(1))
\]
and the Hasse invariant
\[
H^2(G_v, \Qp(1)) \xto{\sim} \Qp
\]
we define
\[
\phi_v:X(F_v) \to \Qp
\]
by
\[
\phi_v(a) = \log \chi \cup j_v(a).
\]
We define a \emph{$p$-adic local N\'eron function} to be a function
\[
X(F_v) \to \Qp
\]
which satisfies axioms analogous to those which define the real N\'eron function (see segment \ref{Ner0} below).
Our main goal in section \ref{S4} is Theorem \ref{Ner0}:
\subsection*{Theorem}
The function $\phi_v$ is a $p$-adic local N\'eron function. 

\medskip
Consider Weierstrass coordinates $x,y$ in which $X$ is given by
\[
y^2 +a_1 xy +a_3y = x^3 +a_2x^2 +a_4x +a_6.
\]
Among the three axioms which define a N\'eron function, verification of the formula
\[
\phi_v(2a) = 4 \phi_v(a) - \log |(2y+a_1x +a_3)(a)|_v
\]
is hardest. This is accomplished via an elaborate computation which takes place on the profinite level and which culminates in the theorem of segment \ref{d18}.

\segment{11_a}{}
Let $\cX=\cE\setminus O$, where $\cE$ is the regular minimal model of an elliptic curve with semi-stable reduction everywhere, let $\al$ be the global $1$-form given by
\[
\al = \frac {dx}{2y +a_1x +a_3}
\]
in Weierstrass coordinates, and let $\be$ be the meromorphic form
\[
\be = x\al.
\]
Let $b$ be the integral tangent vector at $O$ dual to $\al(O)$. Let $S$ denote the set of primes of bad reduction for $\Ee$ and for each $l\in S$, let 
$N_l=\ord_l(\D_{\cE})$, where $\D_{\cE}$ is the minimal discriminant.
Define a set
\[
W_l:=\set {(n(N_l-n)/2N_l)\log l }  { 0\leq n< N_l },
\]
and for each $w=(w_l)_{l\in S} \in W:=\prod_{l\in S} W_l$, define
\[
\|w\|=\sum_{l\in S} w_l.
\]
Our main result in section \ref{S6} is as follows.
\subsection*{Theorem}
Suppose $\cE$ has rank zero and that $\Sha_E[p^{\infty}]<\infty$, and let $p$ be an odd prime of good reduction. 
With assumptions as above
\[
\cX(\Z_p)_2=\bigcup_{w\in W} \Psi(w),
\]
where
\[
\Psi(w):=\set{ z\in \cX(\Z_p) }{ \log (z)=0, \ D_2(z)=\|w\| }.
\]
Here,
\[
\log (z)=\int_b^z \a
\]
and
\[
D_2(z)=\int_b^z \al \be,
\]
are Coleman (iterated) integral functions on $\cX(\Q_p)$.

\segment{}{}
Let us sketch the proof of theorem \ref{11_a}.
It follows from theorem 5.2 of Silverman \cite{SilvermanComputing} that the local height at primes $v \neq p$ takes values in the finite set $W_v$; by theorem \ref{Ner0}, this applies to the image of $\Xx(\ZZ_v)$ under
\[
\Xx(\ZZ_v) \to H^1(G_v, \Qp(1)) \xto{\phi_v} \Qp
\]
where $\phi_v$ now denotes the map
\[
c \mapsto \log \chi \cup c.
\]
As we explain in segment \ref{11_b}, the map $\loc_p$ at level $2$ factors as
\[
\Sel^2(\Xx) \to H^1_f(G_p, \Qp(1)) \inj H^1_f(G_p, U_2);
\]
it is here that we use the assumption about the rank. Drawing on global reciprocity, we find that the image of $\Sel^2(\Xx)$ in $H^1_f(G_p, \Qp(1))$ is given by
\[
\set{\eta  }{ \phi_p(\eta) = \| w\| \mbox{ for some } w\in W}.
\]
As we explain in remark \ref{11_d}, this hints at the possibility of a certain nonabelian reciprocity law, an idea carried further by the third author in \cite{KimRecip}. This also translates into a proof of the theorem, through the unipotent Bloch-Kato exponential.

Armed with theorem \ref{11_a} we are able to verify conjecture \ref{13i} for the prime $p=5$ for 256 semi-stable elliptic curves of rank zero from Cremona's table. We also extend our discussion of punctured elliptic curves with a brief treatment of the rank-one case.

\segment{}{}
In section \ref{S5} we consider the thrice punctured line over $\ZZ$. Of course, in this case the set of $\ZZ$-points is empty, so the conjecture holds at level $n$ when $\Xx(\Zp)_n = \emptyset$. Our results may be summarized as follows.
\subsection*{Proposition}
Let $\cX=\P^1\setminus \{0,1,\infty\}$.  Then
$$\cX(\Z_p)_1=\phi$$
if $p\equiv 2 \mod 3$. If $p\equiv 1 \mod 3$, then
$$\cX(\Z_p)_2=\phi$$
if the value of the $p$-adic dilogarithm $Li_2(z)$ at a sixth root of 1 is non-zero. 

\medskip\noindent
We also report on computations showing that indeed $\Li_2(\ze_6) \neq 0 $ in the range
\[
3 \le p \le 10^5.
\]

\segment{}{}
 In section \ref{S7} we discuss curves of genus $\ge 2$. We consider as an example the Fermat curve $X_l$ given by 
 \[
 x^l+y^l = z^l.
 \]
 We find that if the Tate-Shafarevich group of the Jacobian fulfills its conjectured finiteness, if $l=5$ or $7$, and if $p \not \equiv 1 \mod l$, then conjecture \ref{13i} holds at level $1$. We also show how, starting with a punctured elliptic curve which fulfills conjecture \ref{13i} at level $2$ we can construct a curve of higher genus which fulfills the conjecture at level $2$ as well.

\segment{}{}
Finally, in section \ref{S8} we turn to the $S$-integral variant of our conjecture mentioned above. We apply this to the thrice punctured line, concluding that here the conjecture holds for $S = \{2\}$ and $p = 3,5,6$ in depth $2$.

\subsection*{Acknowledgements}
M.K. is grateful to John Coates, Henri Darmon, Kazuya Kato, Florian Pop, and Andrew Wiles for a continuous stream of discussions on the topic of this paper. He is also grateful to Shinichi Mochizuki whose question prompted a precise formulation of the conjecture, and to Yuichiro Hoshi for a kind and detailed reply to a question about a pro-$p$ analogue. We would like to thank the referee for many helpful comments.

\section{Selmer schemes with stringent local conditions}
\label{S2}  

\segment{13a}{}
We let $\Xx \to \Spec \ZZ$ denote a \textit{regular minimal $\ZZ$-model of a hyperbolic curve over $\QQ$}. By this we mean one of the following.
\begin{itemize}
\item[--] $\Xx = \PP^1\setminus \Dd$ where $\Dd$ is a reduced horizontal divisor with at least three $\CC$-points. In this case we let $\Xx' = \PP^1$.
\item[--] The regular minimal model of a compact smooth curve of genus $\ge 2$. We let $\Xx' = \Xx$.
\item[--] The complement of a non-empty reduced horizontal divisor $\Dd$ inside a regular minimal model $\Xx'$ of a compact smooth curve of genus $\ge 1$.
\end{itemize}
We fix a ``base-point'' $b$ of $\Xx$. In all three cases $b$ may be a $\ZZ$-valued point. In the first and third cases, suppose that $\Dd \subset \Yy$ with $\Yy \subset \Xx'$ open and $\Yy \to \Spec \ZZ$ smooth, so that in particular, $\Om^1_{\Xx'/\Spec \ZZ}|_\Dd$ is invertible. Then we allow $b$ to be an ``integral tangent vector'', by which we mean a nowhere vanishing section of the tangent sheaf
\[
\Tt_{\Xx'/\Spec \ZZ}|_\Dd = \Om^{1\lor}_{\Xx'/\Spec \ZZ}|_\Dd
\]
to $\Xx'$ along $\Dd$.

\segment{13b}{}
Let $p$ denote an odd
prime of good reduction. We then have the unipotent $p$-adic \'etale fundamental group $U$ of $\Xx_\QQ$ at $b$ constructed by Deligne \cite{Deligne89}.  We denote its descending central series by $U = U^1 \supset U^2 \supset \cdots$, and the associated quotients by $U_n = U/U^{n+1}$. We also have, for every $x\in \Xx(\ZZ)$, the path torsor $P(x)$, and corresponding quotients $P_n(x)$. We let $T$ denote a finite set of primes which contains all primes of bad reduction for $\Xx'$ and for $\Dd$, plus the auxiliary prime $p$. Let $\QQ_T$ denote the extension of $\QQ$ which is maximal for the property of being unramified outside of $T$, and let $G_T$ denote the Galois group of $\QQ_T$ over $\QQ$. Then as explained in \S2 of \textit{Selmer varieties} \cite{kimii}, $U_n$ possesses a $G_T$-action, and $P_n(x)$ bears the structure of a $G_T$-equivariant $U_n$-torsor, with $G_T$ acting as usual on the left, but $U_n$ acting on the right.

\segment{13c}{}
For each prime $v$, we fix an embedding $\QQ_T \subset \overline{\QQ_v}$ in the algebraic closure of $\QQ_v$. This gives us for every $v$ a map
\[
G_v \to G_T
\]
from the total Galois group of $\QQ_v$ (which, for $v \notin T$, factors through $\hat \ZZ$). This also gives us an isomorphism of $U_n$ with the unipotent fundamental group of $\Xx_{\QQ_v}$, which we continue to denote by the same symbol. For $y \in \Xx(\ZZ_v)$, we have the local path torsor $P_n(y)$, a $G_v$-equivariant $U_n$-torsor. For $y \in \Xx(\Zp)$, the associated torsor $P_n(y)$ is moreover \textit{crystalline} in the sense of \S2 of \textit{Selmer varieties}; as explained there, this follows from Olsson \cite{OlssonTowards}.

\segment{13d}{}
For each prime $v$ there is an affine, finite type $\Qp$-scheme $H^1(G_v, U_n)$ parametrizing $G_v$-equivariant $U_n$-torsors. For $v=p$ there's a closed subscheme
\[
H^1_f(G_p, U_n) \subset H^1(G_p, U_n)
\]
which parametrizes those torsors which are crystalline. There is also the global $H^1(G_T, U_n)$, an affine finite type $\Qp$-scheme parametrizing $G_T$-equivariant $U_n$-torsors, and for each $v$, a map of $\Qp$-schemes
\[
\m{loc}_v:H^1(G_T, U_n) \to H^1(G_v, U_n)
\]
in terms of which we define $H^1_f(G_T, U_n)$ to be the preimage $\m{loc}_p\inv(H^1_f)$ of $H^1_f(G_p, U_n)$ under $\m{loc}_p$. These fit into commuting squares like so.\footnote{Technically speaking, while the map $\loc_v$ appearing in the diagram is a morphism of $\Qp$-schemes, the vertical maps $j$, $j_v$ are just maps of sets into the sets of $\Qp$-points of the varieties below.}
\[
\xymatrix{
\Xx(\ZZ) \ar[d]_j \ar[r] & \Xx(\ZZ_v) \ar[d]^{j_v} \\
H^1_f(G_T, U_n) \ar[r]_{\m{loc}_v} & H^1(G_v, U_n)
}
\]
The vertical map $j$ is called the \emph{global unipotent Kummer map}, and its local counterpart $j_p$ is called the \emph{local unipotent Kummer map}. As above, we refer the reader to \S2 of \textit{Selmer varieties} \cite{kimii} for the details of these constructions. 

\segment{0211a}{Proposition}
Let $v$ be a prime $\neq p$. Then the subset $\Im j_v$ of the rational points of $H^1(G_v, U_n)$ is finite.

\begin{proof}
See Kim-Tamagawa
\cite{KimTamagawa}.
\end{proof}

\segment{13e}{Remark}
For $v \notin T$ a prime of good reduction, we have $\Im j_v = 0$; see the proof of Corollary 0.3 in \S2 of loc. cit.

\segment{13f}{Definitions}
We define the \emph{Selmer scheme of $\Xx$} to be the (infinite) intersection
\[
\Sel^n(\Xx) := \bigcap_{v \neq p} \m{loc}_v\inv \big( \Im j_v \big)
\]
with scheme structure defined by the sum of the corresponding ideals. We also refer to $H^1_f(G_p, U_n)$ as the \emph{local Selmer scheme of $\Xx$ near $p$}. As $n$ varies, these form two towers:
\[
\xymatrix{
\vdots \ar[d] & \vdots \ar[d] \\
\Sel^2(\Xx) \ar[d] \ar[r]^-{\m{loc}_p} & H^1_f(G_p, U_2) \ar[d] \\
\Sel^1(\Xx)  \ar[r]^-{\m{loc}_p} & H^1_f(G_p, U)  \\
}
\]
compatible with the maps $\m{loc}_p$ as well as $j$ and $j_p$. Thus, if we set
\[
\Xx(\Zp)_n := j_p\inv \big( \m{loc}_p (\Sel^n(\Xx) \big)
\,,
\]
we obtain a non-increasing sequence of refinements
\[
\Xx(\Zp) \supset \Xx(\Zp)_1 \supset \Xx(\Zp)_2 \supset \cdots \supset \Xx(\ZZ)
\]
of the set of $\Zp$-points, containing the set of global points. We say that $p$-adic points which are contained in $\Xx(\Zp)_n$ are \emph{cohomologically  global of level $n$}, or \emph{weakly global of level $n$}.

\segment{13g}{}
Our first task is to remove the apparent dependence on $T$. 

\subsection*{Lemma}
Let $\G$ and $U$ be topological groups with $\G$ acting continuously on $U$.
Let $N\subset \G$ be a closed normal subgroup. Then there is an exact sequence of pointed sets
$$1\to H^1(\G/N, U^N) \xto{i} H^1(\G,U) \xto{r} H^1(N,U).$$

\begin{proof}
Recall that continous cohomology is defined (\cite{kimi}, section 1) as
$$H^1(\G,U)=U\bs Z^1(\G,U),$$
where
$Z^1(\G,U)$ consists of the continuous maps
$c:\G\rTo U$ such that
$$c(g_1g_2)=c(g_1)g_1c(g_2),$$
while $(uc)(g)=uc(g)g(u^{-1})$ for $u\in U$ and $c\in Z^1(\G,U)$.

It is clear that $r\circ i$ sends everything to the base-point.
Assume $r(c)=0$ for a continuous cocycle $c:\G\rTo U$.
So there is a $u\in U$ such that $c(n)=un(u^{-1})$ for all $n\in N$.
Define
$$b(g)=u^{-1}c(g)g(u),$$
a cocycle in the same $U$-orbit as $c$. Then
$$b(n)=u^{-1}c(n)n(u)=u^{-1}un(u^{-1})n(u)=e$$
for all $n\in N$. Thus,
$$b(gn)=b(g)gb(n)=b(g)g(e)=b(g)$$
for all $g\in \G$ and $n\in N$. Since this also implies $b(ng)=b(gg^{-1}ng)=b(g)$, we get$$nb(g)=b(n)nb(g)=b(ng)=b(g)$$
for all $g\in \G$ and $n\in N$. That is, $b$ factors to a cocycle
$$\bar{b}:\G/N \rTo U^N,$$
which is  continuous since $\G/N$ has the quotient topology.
\end{proof}

\subsection*{Proposition}
If
$T'$ and $T$ are two finite sets of primes that contain all primes of bad reduction and $p$, the natural restriction maps
$$H^1_f(G_T, U_n)\inj H^1_f(G_{T\cup T'}, U_n)\hookleftarrow H^1_f(G_{T'}, U_n)$$
induce isomorphisms of Selmer schemes. 

\begin{proof}
We need only consider an enlargement of $T$ to $T'\supset T$. We work with points with values in an arbitrary $\Qp$-algebra, which we will omit from the notation.
We will provisionally put the sets of primes into the notation, as in
$\Sel^n_T(\Xx)$.
Clearly $\Sel^n_T(\Xx) \inj \Sel^n_{T'}(\Xx)$.
Recall that $T$ contains already all primes of bad reduction and $p$. In particular, the action of $G_v$ for every prime $v\in T'\setminus T$  on $U_n$ is unramified.
Thus, the image of $\Xx(\ZZ_v)$ in $H^1(G_v, U_n)$ is trivial (\S\ref{0211a}).
That is, when $v\in T'\setminus T$, for a cohomology class in $c\in H^1(G_{T'}, U_n)$, the condition of locally belonging to the image of $j_v$ is actually the same as triviality at $v$. Thus, $c$ goes to zero under any of the restriction maps
$$H^1(G_{T'}, U_n) \to H^1(I_v, U_n).$$
Since the $I_v$ act trivially on $U_n$, this implies that $c$ goes to zero under the restriction map
$$H^1(G_{T'}, U_n) \to H^1(N, U_n),$$
where $N\subset G_{T'}$ is the subgroup generated by $I_v$ for $v\in T'\setminus T$.
According to lemma \ref{13g}, it follows that $c$ comes from $H^1(G_T, U_n)$. By the commutativity of the triangle
\[
\xymatrix{
H^1(G_T, U_n) \ar[rr] \ar[dr]_-{\m{loc}_v} && H^1(G_{T'}, U_n) \ar[dl]^-{\m{loc}_v} \\
& H^1(G_v, U_n)
}
\]
the local conditions remain the same for both spaces, and hence,
\[
\Sel^n_T(\Xx) \simeq \Sel^n_{T'}(\Xx).
\qedhere
\]
\end{proof}

\subsection*{Corollary}
The subset  $\Xx(\Zp)_n\subset \Xx(\Zp)$ is independent of the choice of the set of primes $T$.

\segment{13h}{}
Now we consider the possibility of a change of base-point from $b$ to $c$.
For this discussion, we will write
$U(b)$ and $U(c)$ for the prounipotent $p$-adic \'etale
fundamental groups with base-points at $b$ and $c$ respectively. Denote by $P(b,x)$ the torsor of prounipotent $p$-adic \'etale paths from $b$ to $x$ ($P(x)$ above). Now, given
 any torsor $W$ for $U(b)$, we get the torsor
\[
W^c:=W\times_{U(b)} P(b,c)=[W\times P(b,c)]/U(b).
\]
Here the action of $u\in U(b)$ takes
$(w, \g)\in W\times P(b,c)$ to $(wu, u^{-1}\g)$. This construction defines  a map
from the groupoid of $U(b)$ torsors to the groupoid of $U(c)$-torsors.

\subsection*{Lemma}
If $b$ and $c$ are both integral, then 
\[
W\mapsto W^c
\]
maps unramified torsors at $ v\notin T$ to unramified torsors, and  crystalline torsors  at $p$ to crystalline torsors. 

\begin{proof}
The condition of being unramified at $v$ is given by triviality under the restriction
map
\[
H^1(G_v, U) \to H^1(I_v, U),
\]
while the crystalline condition is given by triviality under the  map
\[
H^1(G_p, U) \to H^1(G_p, U(B_{cr})).
\]
But since $P(b,c)$ is itself unramified at $v\notin T$ and crystalline at $p$,
both conditions are preserved by the functor.
\end{proof}

\noindent
That is, we are assured of an isomorphism
\[
(\cdot)^c: H^1_f(G_T, U(b) )\simeq H^1_f(G_T, U(c)).
\]
Meanwhile, since
\[
(P(b,x))^c=P(c,x),
\]
torsors of paths are preserved under the functor. So we conclude

\subsection*{Proposition}
The functor $W \mapsto W^c$ induces isomorphisms of local and global Selmer schemes commuting with the corresponding localization maps $\m{loc}_p$ and Kummer maps $j$ and $j_p$.

\subsection*{Corollary}
The subset  $\Xx(\Zp)_n\subset \Xx(\Zp)$ is independent of the choice of base-point $b$.

\section{The Conjecture and its context}
\label{S3}

\segment{13i}{}
We preserve the situation and notation of \S1. In particular, $\Xx$ denotes a \textit{minimal $\ZZ$-model of a hyperbolic curve over $\QQ$} as in Segment \ref{13a}. In his lectures at the IH\'ES in February of 2012, M.K. proposed the following.

\subsection*{Conjecture}
Equality $\Xx(\Zp)_n = \Xx(\ZZ)$ is obtained for large $n$.

\segment{v2a}{Remark}
Although we would expect a suitable generalization of our conjecture to hold over general number fields, the exact statement is not entirely clear, and we do not go into this issue in this paper. See Dan-Cohen \cite{mtmueII} for the case of the thrice punctured line.

\segment{13j}{}
Recall that $j$, $j_p$ denote the global and local Kummer maps, respectively (\ref{13d}). Alongside conjecture \ref{13i}, we consider the following statements.
\begin{enumerate}
\item[(SGK)] \textit{Surjectivity of the global Kummer map.} The global Kummer map $j$ defines a surjection 
\[
X(\ZZ) \surj \set{P \in \Sel^n(\Xx)}{\m{loc}_p(P) \in \Im j_p}
\]
onto the set of torsors which are \textit{geometric everywhere locally}, for large $n$.
\item[(ILK)] \textit{Injectivity of the local Kummer map.} Suppose $x \in \Xx(\ZZ)$ and $y \in \Xx(\Zp)$. If $j_p(x) = j_p(y)$ for all $n$ then $x=y$.
\end{enumerate}
Trivially, we have the implications
\[
\mbox{SGK + ILK } \Rightarrow \mbox{ Conjecture \ref{13i} } \Rightarrow \mbox{ ILK.}
\]

\segment{13k}{Relationship to Tate--Shafarevich and Section Conjectures}
We now discuss the relationship between \ref{13j}(SGK), finiteness of Sha, and the Grothendieck section conjecture. Let $X$ be a proper hyperbolic curve over $\QQ$, let $b \in X(\QQ)$, fix an algebraic closure $\bar \QQ$ of $\QQ$, and let $G_{\QQ}$ denote the Galois group of $\bar \QQ / \QQ$. For each $x \in X(\QQ)$ we let $\hat P(x)$ denote the $\pi_1^\et(X_{\bar \QQ},b)$-torsor associated to $x$. Recall that the Grothendieck section conjecture states that
\[
\hat j: x \mapsto \hat P(x)
\]
defines a bijection
\[
X(\QQ) = H^1\big(G_\QQ, \pi_1^\et(X_{\bar \QQ}, b)\big).
\]
The surjectivity of $\hat j$ bears an obvious relationship to statement \ref{13j}(SGK). When we replace $\pi_1^\et(X_{\bar \QQ},b)$ by its prounipotent completion $U$, the cohomology set becomes a positive dimensional variety, so surjectivity ceases to be plausible; $j$ may nevertheless surject onto those cohomology classes which are everywhere locally geometric. This is motivated in part by the case of elliptic curves and the conjectured finiteness of $\Sha$, through the following basic proposition.

\ssegment{14a}{}
Let $E$ be an elliptic curve over $\QQ$. As above, we fix a decomposition group $G_v \subset G_\QQ$ at every prime $v$, giving rise to a localization map
\[
\loc_v: H^1\big(G_\QQ, H^\et_1(E_{\bar \QQ}, \Qp)\big) \to H^1\big(G_v, H_1^\et(E_{\bar \QQ}, \Qp)\big)
\,.
\]
Let $j^\Qp$, $j_v^\Qp$ denote the global and local $\Qp$-linearized Kummer maps, as in the following square.
\[
\xymatrix{
E(\QQ)_{/p}\otimes_\Zp \Qp \ar[r] \ar[d]_{j^\Qp} & E(\QQ_v)_{/p}\otimes_\Zp \Qp \ar[d]^{j^\Qp_v} \\
H^1\big(G_\QQ, H_1(E_{\bar\QQ}, \Qp)\big) \ar[r] & H^1\big(G_v, H_1(E_{\bar\QQ}, \Qp)\big) 
}
\]
Here the subscript $/p$ denotes $p$-adic completion. 

\subsection*{Proposition}
Suppose the $p$-part of $\Sha(E)$ is finite. Then the global (abelian) $\Qp$-linearized Kummer map $j^\Qp$ defines a bijection
\[
E(\QQ) \otimes_{\Zp} \Qp =
 \set{P\in H^1\big(G_\QQ, H^\et_1(E_{\bar \QQ}, \Qp)\big)  }{\loc_v(P) \in \Im j_v^\Qp \mbox { for all primes }v}
\]
between the vector space of linear combinations of rational points and the classical $p$-adic Selmer group.

\begin{proof}
The finiteness of the $p$-part of $\Sha$ implies that the product of the localization maps induces an injection
\[
\tag{*}
\varprojlim H^1(G_\QQ, E)[p^n] \otimes \Qp
 \inj
  \prod_v \varprojlim H^1(G_v, E)[p^n] \otimes \Qp.
\]
Recall that there's a Galois-equivariant isomorphism
\[
\Qp \otimes \varprojlim E[p^n](\Qbar) = H_1^\et(E_{\Qbar}, \Qp).
\]
We consider the inverse system of short exact sequences
\[
\xymatrix
{
0 \ar[r] & E[p^n] \ar[r] & E \ar[r]^{p^n} & E \ar[r] & 0 \\
0 \ar[r] & E[p^{n+1}] \ar[u] \ar[r] & E \ar[u]_p \ar[r]_{p^{n+1}} & E \ar@{=}[u] \ar[r] & 0.
}
\]
Taking $\Qbar$-valued points followed by invariants by $G_\QQ$, we obtain an inverse system of short exact sequences
\begin{small}
\[
\xymatrix
{
0 \ar[r] & 
E(\QQ)/p^n \ar[r] & 
H^1 \big( G_\QQ, E[p^n](\Qbar) \big) \ar[r] & 
H^1 \big( G_\QQ, E(\Qbar) \big) [p^n] \ar[r] & 
0 
\\
0 \ar[r] & 
E(\QQ)/p^{n+1} \ar[r] \ar[u] & 
H^1 \big( G_\QQ, E[p^{n+1}](\Qbar) \big) \ar[r] \ar[u] & 
H^1 \big( G_\QQ, E(\Qbar) \big) [p^{n+1}] \ar[r] \ar[u] & 
0. }
\]
\end{small}
Taking inverse limits and tensoring with $\Qp$, we obtain the top row in the following diagram:
\begin{small}
\[
\xymatrix
@ C=3ex
{
0 \ar[r] &
 E(\QQ)_{/p}\otimes \Qp \ar[r]^-{j^\Qp} \ar[d] & 
 H^1\big(G_{\QQ}, H^\et_1(E_{\Qbar}, \Qp) \big) \ar[r] \ar[d] & 
 \varprojlim
  H^1\big( G_\QQ, E(\Qbar) \big)[p^n] \otimes \Qp
  \ar[r] \ar[d] & 
0 \\
0 \ar[r] & 
E(\QQ_v)_{/p}\otimes \Qp \ar[r]_-{j_v^\Qp} & 
H^1\big(G_{\QQ_v}, H^\et_1(E_{\overline \QQ_v}, \Qp) \big) \ar[r] & 
\varprojlim
 H^1 \big( G_{\QQ_v}, E(\overline \QQ_v) \big)[p^n] \otimes \Qp
  \ar[r] & 
0 ;
}
\]
\end{small}
\noindent
repeating the procedure with $\QQ_v$ in place of $\QQ$ gives us the rest of the diagram. Varying the place $v$ and using the injectivity (*),  we obtain an exact sequence like so
\[
0 \to 
 E(\QQ)_{/p}\otimes \Qp 
  \to 
 H^1 \big(G_{\QQ}, H^\et_1(E_{\Qbar}, \Qp) \big)
  \to
 \prod_v 
 \frac
 {H^1\big(G_{\QQ_v}, H^\et_1(E_{\Qbar}, \Qp) \big)}
 { \big(E(\QQ_v)_{/p}\otimes \Qp \big)} 
 \;.
\]
By the Mordell-Weil theorem, we have 
\[
E(\QQ)_{/p}\otimes \Qp = E(\QQ) \otimes \Qp,  
\]
so the proposition follows.
\end{proof}

\subsection*{Corollary}
We have $ E(\QQ)\otimes \Qp = \Sel^1(E) $. In particular, our $\Sel^1(E)$ is equal to the classical $p$-adic Selmer group.

\begin{proof}
For $v \neq p$, the geometricity condition
\[
\loc_v(P) \in \Im j_v
\]
is actually equivalent to the (a priori weaker) condition
\[
\loc_v(P) \in \Im j_v^\Qp
\]
since
\[
E(\QQ_v)_{/p} \otimes \Qp = 0.
\]
On the other hand at $v=p$ the condition $\loc_p(P) \in \Im j_p^\Qp$ is equivalent to $\loc_p(P)$ being crystalline according to Example 3.11 of Bloch-Kato \cite{BlochKato}.
\end{proof}

\section{The unipotent Albanese map and local height on elliptic curves}
\label{S4}

\segment{}{Setup and statement}

\ssegment{d2}{}
As explained in the introduction, our goal here is to investigate  a relation between local heights and Albanese maps,
with a view towards applying it to the computation of some simple Selmer schemes. The relation over a finite extension of $\Q_p$ was noticed earlier following the paper \cite{BalakAppendix} by its authors and
Amnon Besser \cite{BalBesColGr}. The main purpose here will be to work out a precise relation over $\Q_l$ for $l\neq p$, when the curve has bad reduction.

\ssegment{d3}{}
Fix an odd prime $p$, let $F_v$ be an unramified finite extension of $\QQ_l$ for $l\neq p$, and let $(E_v,O)$ be an elliptic curve over $F_v$ written
in Weierstrass minimal form
\[
Z_0Z_2^2+a_1Z_0Z_1Z_2+a_3Z_0^2Z_2=Z_1^3+a_2Z_0Z_1^2+a_4Z_0^2Z_1+a_6Z_0^3.
\]
Let $X=E_v\setminus \{O\}$  with
 equation
\[
y^2+a_1xy+a_3y=x^3+a_2x^2+a_4x+a_6.
\]
We let $b$ be the tangent vector to $E$ at $O$ dual to the invariant differential form
\[
dx/(2y+a_1x+a_3),
\]
which we will use as the main base-point for fundamental groups. Let $z=(-x/y)$, which is  a $b$-{\em compatible} uniformizing element at $O$ in that $(d/dz)|_O=b$. (We refer to Silverman \cite{SilArEl}, chapter 4, for this and other assertions about the coordinates on the Weierstrass minimal model.)

\ssegment{d20}{}
Let $\log$ denote the $p$-adic logarithm normalized so that $\log(p) = 0$. The $p$-adic logarithm
\[
\log \chi: G_v^\m{ab} \to \Qp
\]
of the $p$-adic cyclotomic character may be regarded as an element of $H^1(G_v, \Qp)$. Recall that the cup product defines a $\Qp$-valued pairing
\[
H^1(G_v, \Qp) \times H^1(G_v, \Qp(1)) \to H^2(G_v, \Qp(1)) \xto \cong \Qp.
\]
Let $\rec$ denote the reciprocity map of abelian class field theory
\[
F_v^* \to G_v^\m{ab}
\]
and recall that $l$ denotes the residue characteristic of $F_v$. Then for $a \in F_v^*$ we have the formula
\[
\log \chi \cup k(a) = \log(\chi(\rec_v(a))) = \log l^{v(a)} = -\log |a|_v.
\]

\ssegment{d21}{Proposition}
We have
\[
H^0(G_v, U_1) = H^1(G_v, U_1) = H^2(G_v, U_1) = 0.
\]

\begin{proof}
We start with $H^0$: by the weight--monodromy theorem, proved for abelian varieties by Grothendieck in \cite[Expos\'e IX]{SGAIV}, the inertia fixed part of $E[p^n]$ has Frobeinus weight $-2$, so in particular has no Frobenius-fixed part, whence the vanishing. The vanishing of $H^2$ then follows by local Tate duality \cite{NeukirchCoh}, since $V_p(E_v)$ is self-dual. Since the $v$-adic absolute value of $p^n$ is $1$, it follows from \cite[Theorem 7.3.1]{NeukirchCoh} that the Euler characteristic is zero; combined with the vanishing of $H^0$ and $H^2$, this implies the vanishing of $H^1$.
\end{proof}

\ssegment{bgu2ta}{}
Recall that the $G_v$-equivariant extension
\[
1 \to \Qp(1) \to U_2 \to U_1 \to 1
\] 
gives rise to an exact sequence of pointed sets
\[
H^0(U_1) \to H^1(\Qp(1)) \xto \al H^1(U_2) \to H^1(U_1).
\]
The vanishing of the extreme terms implies that $\al$ is bijective, so that we can choose a cocycle representing $j_v(x)$ which takes values in $\Qp(1)$. Thus, we get a function
\[
\phi_v: X(F_v)\rTo \Q_p
\]
via the formula
\[
\phi_v(a)=\log\chi \cup j_v(a).
\]

\ssegment{Ner0}{}
 We define a \emph{$p$-adic local N\'eron function} to be a function 
\[
\la:E_v(F_v) \to \Qp
\]
which satisfies the following properties with respect to the coordinates $x$, $y$.
\begin{enumerate}
\item[(i)] $\la$ is continuous on $E_v(F_v) \setminus \{O\} $ and bounded on the complement of any $v$-adic neighborhood of $O$.
\item[(ii)] The limit
\[
\underset{a \to 0}\lim \big( \la(a) - \frac{1}{2}\log|x(a)|_v \big)
\]
exists.
\item[(iii)] For all $a \in E_v(F_v)$ with $[2]a \neq 0$, 
\[
\la([2]a) = 4\la(a) - \log|(2y+a_1x+a_3)(a)|_v.
\]
\end{enumerate}

\subsection*{Theorem}
The function $\phi_v$ is a $p$-adic local N\'eron function.

\bigskip

The proof of theorem \ref{Ner0} appears in segment \ref{phiv2x} below.

\segment{pi2}{Construction of $\pi_{[2]}$-tower}

\ssegment{d4}{}
We write here $\pip$ for
\[
\pi^{(p)}(\bX,b)/[ \pi^{(p)}(\bX,b) , [  \pi^{(p)}(\bX,b),   \pi^{(p)}(\bX,b)] ],
\]
the quotient of the pro-$p$ fundamental group of $\bX=X\otimes \bF_v$
by the third level of its lower central series. We will need to consider   different base-points $w$ below, in which case we denote the group by
$\pip(w)$. Similarly, the pushout to $\pip(w)$ of the homotopy class of maps from $w$ to $y$ will be denoted by $\pip(w,y)$:
\[
\pip(w,y):=\pi^{(p)}_1(\bX; w, y)\times_{\pi^{(p)}_1(\bX, w)}\pip(w).
\]
Note that when $b$ is replaced by $\lambda b$ for $\lambda \in F_v$, then the compatible uniformizer is changed to $z/\lambda$.

We note that $\pip$ fits into an exact sequence
\[
0 \to  
Z    \to \pip \to T_pE \to 0
\]
where
\[
Z\simeq \Z_p(1)
\]
is generated by $[e,f]$ for any lift $\{e, f\}$ of a basis for $T_pE$. As in Lemma 1.1 of \cite{KimMassey}, this exact sequence has a Galois-equivariant splitting which extends also to a splitting of the sequence
\[
1 \to \Qp(1) \to U_2 \to V_pE \to 0.
\]
Thus, as in \cite[p. 730]{KimMassey}, we will write a cocycle
\[
c: G_v \to U_2
\]
as
\[
c=c_2c_1,
\]
where $c_2$ takes values in $\Qp(1)$, $c_1$ is a cocycle with values in $V_pE$, and
\[
dc_2 = -(1/2)c_1 \cup c_1.
\]

We wish to compute the group $\pip$ using theta groups. The result is stated in proposition \ref{d12} below.

\ssegment{d5}{}
Let
$
D_0 := [p^n]^*[O],
$
the sum of all points of $E_v[p^n]$. We write $\sim$ for linear equivalence of divisors. We claim that 
\[
D_0 \sim p^{2n}[O].
\]
To see this we base change to an algebraically closed field, write 
\[
D_0=\sum_{j=0}^{p^{2n}-1} [z_j],
\]
and remember the isomorphism of group schemes 
\[
E_v \xto{\sim} \Pic^0 E_v
\]
\[
z \mapsto [z] -[O],
\]
from which
\begin{align*}
D_0-p^{2n}[O] 
	&= \Big( \sum_j [z_j] \Big) -p^{2n}[O] \\
	&= \sum_j \big([z_j] -[O] \big) \\
	&\sim \Big[ \sum_j z_j \Big] -[O] \\
	&= [O] -[O] \\
	&= 0.
\end{align*}

Let $\cH_n :=\O(p^n[O])$. Then we have
\[
\cH_n^{p^n}
\simeq 
\Oo\Big(p^n\big(p^n[O] \big) \Big)
 = \Oo\big(p^{2n}[O] \big)
  \simeq 
  \O(D_0),
\]
via an isomorphism well-defined up to a constant. 

\ssegment{d6}{}
In general, an isomorphism
\[
\O(A)\simeq \O(B)
\]
must be defined by a rational function $f$ such that
$(f)=B-A$, which takes a section $s\in \O(A)$ and multiplies it by $f$.
We will denote this isomorphism also by $f$:
\[
\O(A)\stackrel{f}{\simeq} \O(B).
\]
When a tangential base-point $w$ at $O$ has been chosen we  normalize all such isomorphisms as follows. Choose a local coordinate $t$ at $O$ so that $(d/dt)|_O=w$. Then we normalize $f$ so that
\[
t^{-\ord_O(f)}f(O)=1.
\]
When this normalization has been fixed, we will refer to the function or the
isomorphism as {\em based} at $w$.
This way, when
\[
\O(A)\stackrel{f}{\simeq} \O(B), \ \ \ \O(B)\stackrel{g}{\simeq} \O(C)
\]
and
\[
\O(A)\stackrel{h}{\simeq} \O(C),
\]
are all based at $w$, then
we can be sure that $h=gf.$ We will be able to deduce the commutativity of various diagrams using this fact. The based function giving the isomorphism
\[
\cH_n^{p^n}\simeq \O(D_0)
\]
given a base-point $w$ will be denoted by $f_w$. More generally, given any function  $g$ such that
\[
(g)=A-B
\]
we will write $g_w$ for the constant multiple of $g$ that is based at the tangent vector $w$.

\ssegment{d7}{}
For our choice of tangent vector $b$, an elementary computation shows that the function $y$ is based, that is, $y\sim z^3$. In the case of the function $f_b\in \Q[x,y]$, clearly, there is a constant multiple $f^{\Z}\in \Z[x,y]$ with the property that $(f^{\Z})_{\infty}=(p^{2n}-1)[O]$
 on the Weierstrass minimal model. But then, since $(z)=(-x/y)=[O]+D$ on the minimal model with $D$ disjoint from $[O]$, we see by comparing divisors that $z^{1-p^{2n}}f^{\ZZ}$ is a unit $h$ in a neighborhood of the section $O$. Thus, its value on $O$ is a unit $u \in \Oo_v^*$. Hence we see that
\[
f_b = u\inv f^\ZZ = u\inv h z^{p^{2n}-1},
\]
with the second equality holding in a neighborhood of $O$. In particular, the formal power series expansion of $f_b$ in the parameter $z$ has coefficients in $\Oo_v$.

\ssegment{d10}{}
Define the subscheme
\[
X_n'\subset \cH_n
\]
as the inverse image of the section $1 \in \G(\O(D_0))$ under the map
\[
\cH_n \xto{(\cdot)^{p^n}} \cH_n^{p^n} \simeq \O(D_0),
\]
where the second isomorphism is given by the function $f_{b/p^n}$. Standard Kummer theory implies that
\[
X_n' \to E_v
\]
is a finite $\mu_{p^n}$ cover (totally) ramified only over $D_0$.
In particular,
\[
r_n: X_n' \to E_v\xto{p^n} E_v,
\]
is a finite cover.

\ssegment{d11}{}
Since the cover $X'_n$ is constructed locally as $\O_E[(f_{b/p^n})^{1/p^n}]$ and 
\[
(u^{p^{2n}}f_{b/p^n})=u+c_2u^2+c_3u^3+\cdots
\]
formally with respect to the uniformizer $u=p^nz$,
 we see that the tangent vector $b/p^n=(d/du)|_O$ lifts to a tangent vector
$b'$ to $X'_n$ at the unique point above $O$. That is, $r_n:X_n \rTo X$ is equipped with an $F_v$-rational
lift of the tangential base-point $b$.
Therefore, there is a $G_v$-equivariant surjective homomorphism
\[
\pip \rTo (X_n)_{b}
\]
that sends the identity to the base-point lift $b'$. Here the subscript $(\cdot)_b$ refers to the tangential fiber functor of Deligne \cite[\S15]{Deligne89}. We will use theta groups to show that this map induces a bijection
\[
\pi_{[2]} \xto{\sim} \varprojlim (X_n)_b.
\]

\ssegment{d8}{}
For each $x\in E_v$, we let $\tau_x:E_v\rTo E_v$ be the translation operator
$\t_x(y)=y+x.$ Recall from section 23 of Mumford \cite{MumfordAbVars} that the theta group 
\[
\cG(\cH_n)
\]
associated to $\cH_n$ is the group scheme over $F_v$ whose $R$-points, for $R$ an $F_v$-algebra, are commuting squares
\[
\xymatrix{
\cH_{n,R} \ar[r]^-g_-\cong \ar[d] &
\cH_{n,R} \ar[d] \\
E_{v,R} \ar[r]_-{\t_x}^-\cong &
E_{v,R}
}
\]
for $x \in E_v(R)$. Since $x$ is determined by $g$, we denote such a square simply by $g$. If $x \in E_v(R)$ then the ideal defining the associated closed subscheme of $E_{v,R}$ is locally principal; we denote the associated Cartier divisor by $[x]$. In this notation, we have
\[
\t^*_{-x}\cH_n \cong \Oo(p^n[x])
\]
isomorphic to $\Oo(p^n[O])$ if and only if $x\in E_v[p^n](R)$. The theta group therefore fits into a short exact sequence
\[
0 \to \Gm \to \Gg(\Hh_n) \xto{\rho} E_v[p^n] \to 0.
\]
Moreover, after forgetting the group structures, the projection $\rho$ admits a section (see segment \ref{d9} below). 

We will also consider a point $g \in \Gg(\Hh_n)$ as an isomorphism
\[
g: \Hh_n \xto{\sim} \tau_{\rho (g)}^* \Hh_n,
\]
in which case composition in $\Gg(\Hh_n)$ is given by the formula 
\[
\tag{*}
g \cdot h = \tau_{\rho(h)}^*(g) \circ h.
\]

\ssegment{3_a}{}
Taking tensor powers defines a map of exact sequences (solid arrow diagram below)
\[
\xymatrix{
0 \ar[r] &
\Gm \ar[r] &
\Gg(\Hh_n^{p^n}) \ar[r] &
E_v[p^{2n}] \ar[r] &
0
\\
0 \ar[r] &
\Gm \ar[u]^-{p^n} \ar[r] &
\Gg(\Hh_n) 
\ar[r]_-{\rho}
 \ar[u]_-{(\cdot)^{\otimes {p^n}}}
 &
E_v[p^n] \ar[r] \ar@{}[u]|\bigcup \ar@/_8pt/@{.>}[ul]^-\be &
0.
}
\]
We now construct a diagonal homomorphism $\be$, as shown, which will make the triangle to its upper right commute. 
\label{3_b}{}
Given $x \in E_v[p^n]$ we let $\psi_x$ denote the canonical isomorphism
\[
\Oo(D_0) \xto{\sim} \tau_x^* \Oo(D_0).
\]
(which is not necessarily compatible with the base point). We note that
\[
\tag{*}
\psi_O = Id_{\Oo(D_0)},
\]
and that if $y \in E_v[p^n]$ is a second point, then the square
\[
\tag{**}
\xymatrix{
\Oo(D_0) \ar[r]^-{\psi_y} \ar[drr]_-{\psi_{x+y}} &
\tau_y^*\Oo(D_0) \ar[r]^-{\tau_y^*(\psi_x)} &
\tau_y^*\tau_x^* \Oo(D_0)  \ar@{=}[d]
\\
&& \tau_{x+y}^* \Oo(D_o)
}
\]
commutes. We define $\be$ by
\[
\be(x) := \tau_x^*f \circ \psi_x \circ f\inv.
\]
Properties \ref{d8}(*), \ref{3_a}(*), and \ref{3_a}(**) combine to show that $\be$ is a homomorphism:
\[
\be(O) = f \circ  \Id \circ f\inv = \Id,
\]
and
\begin{align*}
\be(x)\cdot \be(y)
	&= \tau_y^*(\tau_x^*f \circ \psi_x \circ f\inv )
		\circ \tau_y^*f \circ \psi_y \circ f\inv \\
	&= \tau^*_{x+y} f \circ \tau_y^*(\psi_x) \circ \tau^*_yf\inv
		\circ \tau_y^*f \circ \psi_y \circ f\inv \\
	&= \tau^*_{x+y} f \circ \tau_y^*(\psi_x) \circ 
		 \psi_y \circ f\inv \\
	&= \tau^*_{x+y}f \circ \psi_{x+y} \circ f\inv \\
	&= \be(x+y).
\end{align*}

\ssegment{GnExt}{}
We define
\[
\Gg_n := \big[ (\,\cdot\,)^{\otimes p^n} \big]
\inv \big(\Im \be \big).
\]
This subgroup of the theta group fits into an exact sequence
\[
0 \to \mu_{p^n} \to \Gg_n \xto \rho E_v[p^n] \to 0.
\tag{$*$}
\]
As is customary when doing Galois-theoretic computations, we will often identify $\Gg_n$ with $\Gg_n(\bar F_v)$. 

\subsection*{Lemma}
The finite \'etale covering
\[
r_n: X_n \to X
\]
is Galois with Galois group $\Gg_n$.

\begin{proof}
We claim that the square
\[
\xymatrix{
\Gg(\Hh_n^{p^n})
\\
\Gg(\Hh_n) \ar[u]^{\otimes{p^n}} & E_v[p^n] \ar[ul]_\be 
\\
\Gg_n \ar@{}[u]|\bigcup \ar[ur]_{\rho}
}
\]
commutes. This is an elementary computation which we carry out anyway. We put ourselves in the general setting of a diagram
\[
\xymatrix{
G \ar[r]^\rho  & H \\
G' \ar[u]^{t} \ar[r]_{\rho'} & H' \ar[ul]^\be \ar@{^(->}[u]_u
}
\]
in which the outer square and the upper right triangle commute, and $u$ is injective as shown. If $g \in t\inv(\Im \be)$ then there's a $g' \in G'$ such that
\[
t(g) = \be(\rho'(g')).
\]
We then have
\[
u \rho' g' = \rho \be \rho' g' = \rho t g = u \rho'g,
\]
from which
\[
\rho' g' = \rho' g,
\]
so that 
\[
t(g) = \be \rho' g' = \be \rho' g,
\]
as hoped. 

It follows that the diagram
\[
\begin{tikzcd}
\Hh_n \arrow[loop below, distance = 25]
\arrow[r]
&
\Hh^{p^n} \arrow[loop below, distance = 25] 
\\
\Gg_n \arrow[r]
&
E_v[p^n]
\end{tikzcd}
\]
commutes, in the sense that the map of $G$-sets is linear over the map of groups.

Denote by $Y_n$ the image of $E_v \setminus E_v[p^n]$ under the section $1$ viewed as a map of schemes
\[
E_v \to \Oo(D_0).
\]
Since the action of $ E_v[p^n]$ on $\Oo(D_0)$ given by the liftings $\psi_x$ maps a function $h$ (viewed as a section) to $h \circ \tau_x$, $Y_n$ is stable under the action of $E_v[p^n]$ via $\be$. Hence, $\Gg_n$ acts on $X_n$. 

Further,
\[
\Hh_n \to \Hh_n^{p^n} \cong \Oo(D_0)
\]
is a $\mu_{p^n}$-torsor away from the zero section, where this $\mu^{p^n}$ is exactly the subgroup of $\Gg_n$ mapping to $1$ under the map $(\,\cdot\,)^{\otimes p^n}$. Therefore, $X_n$ is a $\mu_{p^n}$-torsor over $Y_n$ and
\[
X_n/ \mu_{p^n} \cong Y_n.
\]
However, $Y_n$ is isomorphic to $E_v \setminus E_v[p^n]$ equivariantly with respect to the action of $E_v[p^n]$. So
\[
X_n/\Gg_n \cong Y_n/E_v[p^n] \cong E_v\setminus O =X,
\]
which completes the proof of the lemma.
\end{proof}

\ssegment{d12}{}
If we denote by $\tbX_{[2]} \ra \bX$ the  quotient of
the universal covering space of $\bX$ corresponding to $\pip$,
there is a surjective homomorphism
\[
\tag{*}
\Aut(\tbX_{[2]}/\bX)\rOnto \cG_n,
\]
simply because $\Gg_n$ is a central extension of $E_v[p^n]$. The significance of the theta-group for us is that the commutator map
\[
[ \cdot, \cdot ]: \cG_n \times \cG_n \rTo \mu_{p^n}
\]
factors to the Weil pairing
\[
\langle \cdot, \cdot \rangle: E_v[p^n]\times E_v[p^n] \rTo \mu_{p^n},
\]
so is in particular surjective (this is a well known fact, explained for instance in \cite[Chapter XI, Proposition 11.20]{Moonen}). Hence, we also have a surjection
\[
[\Aut(\tbX_{[2]}/\bX), \Aut(\tbX_{[2]}/\bX)] \rOnto \mu_{p^n}.
\]

\subsection*{Proposition}
The surjections \ref{d12}(*) induce an isomorphism of profinite groups
\[
\Aut(\tbX_{[2]}/\bX)\simeq \varprojlim \cG_n,
\]
and hence, a bijection of profinite sets
\[
\pip\simeq \varprojlim (X_n)_b.
\]

\begin{proof} 
It suffices to show injectivity.
So let $g\in \Aut(\tbX_{[2]}/\bX)$ be non-trivial. If
$g$ has non-trivial image in $\Aut(\tbX_{[2]}/\bX)^{ab}\simeq T_pE$, then
clearly there is a map 
\[
\Aut(\tbX_{[2]}/\bX)\rTo \cG_n
\]
which does not send it to zero. So assume
\[
g\in [\Aut(\tbX_{[2]}/\bX), \Aut(\tbX_{[2]}/\bX)].
\]
But then, since
\[
[\Aut(\tbX_{[2]}/\bX), \Aut(\tbX_{[2]}/\bX)]\simeq \Z_p(1)
\]
as a topological group,
the family of surjections
\[
[\Aut(\tbX_{[2]}/\bX), \Aut(\tbX_{[2]}/\bX)]\rOnto \mu_{p^n}
\]
must be separating. 

For the statement about $\pip$, recall that the formula
\[
b'l=\phi (b')
\]
for $l\in \pip$ and $\phi \in \Aut(\tbX_{[2]}/\bX)$
defines an anti-isomorphism from $\pip$ to $\Aut(\tbX_{[2]}/\bX)$.
\end{proof}

\segment{jv2x}{Interaction of local Kummer map with multiplication by $2$}
Our goal here is to derive an explicit formula for the change in $j_v(x)$ when we multiply $x$ by $2$. The result is stated in corollary \ref{d19}.

\ssegment{d9}{}
We can construct a canonical section of the surjection $\rho$ (\ref{GnExt}($*$)) as follows. 
Notice that the automorphism
\[
[-1]:E_v\simeq E_v
\]
lifts to an automorphism
\[
[-1]:\cH_n\simeq \cH_n
\]
that sends a section $\phi(y)$ to $\phi(-y)$. This induces
an involution
\[
i:\cG_n\rTo \cG_n
\]
that sends $g $ to $[-1]\circ g\circ [-1].$

\subsection*{Lemma}
If $R$ is an $F_v$-algebra, then any $R$-valued point of $E_v[p^n]$ has an $R$-valued half. 

\begin{proof}
We recall the proof of this well-known fact. There's a short exact sequence of finite \'etale group schemes
\[
0 \to E_v[2] \to [2]\inv E_v[p^n] \to E_v[p^n] \to 0, 
\]
hence an exact sequence of \'etale cohomologies
\[
[2]\inv E_v[p^n](R) \to E_v[p^n](R) \xto{\delta} H^1(\Spec R, E_v[2]).
\]
The boundary map $\delta$ is a map from a $\ZZ/p^n$-module to a $\ZZ/2$-module, hence, under our assumption that $p$ is odd, necessarily zero.
\end{proof}

Given any element $x\in E_v[p^n]$ and a lift $g \in \cG_n$ of
$-x/2$, the element  $$x':=i(g)g^{-1}$$ is independent of the 
lift $g$, and can be characterized as the unique element of $\cG_n$ lying
over $\t_x$ that is anti-symmetric with respect to $i$, in  that $i(g)=(g)^{-1}$.

Thus, we can write an element
$g\in \cG_n$ uniquely in the form
\[
g=g_2g_1
\]
where $g_1$ is the unique lift of $\rho(g)$ satisfying $i(g_1)=g_1^{-1}$.
Sometimes we  abuse notation and write $g_1$ both for this lift and for
$\rho(g)$.

\ssegment{d13}{}
Recall that the nonabelian cohomology set $H^1(G_v, \pip)$ can be constructed as a set of equivalence classes of continuous cocycles $G_v \to \pip$; it also parametrizes the set of isomorphism classes of \emph{$G_v$-equivariant $\pip$-torsors}. Given a point $x \in X(F_v)$ we write $\hat j(x)$ for the associated class 
\[
\hat j(x):=[\pi^{(p)}_1(\bX; b, x)\times_{\pi^{(p)}_1(\bX, b)}
\pip ]\in H^1(G_v, \pip).
\]
If $c^x$ is an associated cocycle, then composing with the anti-homomorphism 
\[
\pip \surj \Gg_n,
\]
we obtain an anti-cocycle
\[
G_v \to \Gg_n.
\]
which we continue to denote by $c^x$. Explicitly, $c^x$ is constructed as follows. We choose a point $y \in X(\bar F_v)$ such that $p^ny=x$, and a point $z\in X_n(\bar F_v)$ lying above $y$. Then for $\ga \in G_v$, $c^x(\ga) \in \Gg_n$ is determined by the formula
\[
\ga(z) = c^x(\ga)(z).
\]
We remark that the anti-cocycle condition is given by
\[
c^x(\g_1\g_2)=(\g_1c^x(\g_2))c^x(\g_1).
\]
Using the section of $\rho$ constructed in segment \ref{d9}, we can canonically decompose $c^x$ as
\[
c^x = c_2^xc_1^x
\]
with $c_2^x$ taking values in $\mu_n$ and $c_1^x$ the anti-symmetric element lifting $\rho(c^x)$.

\ssegment{d14}{}
We now fix several isomorphisms of line bundles relating to multiplication by $2$ and by $p^n$. There are isomorphisms 
\[
[2]^*(\O([O]))\simeq \O([O]+[x_1]+[x_2]+[x_3])\simeq \O(4[O]);
\]
the first isomorphism is canonical, since
\[
[2]^*(O)=[O]+[x_1]+[x_2]+[x_3]
\]
while for the second isomorphism we may take the one induced by the function
\[
h_b=\frac{2y+a_1x+a_3}{2}.
\]
That is, this function has divisor $3[O]-([x_1]+[x_2]+[x_3])$ and is compatible with the tangent vector $b$. 

There is an isomorphism  
\[
[2]^*(\cH_n)=[2]^*\O(p^n[O])\simeq \O(p^n([O]+[x_1]+[x_2]+[x_3]))
\]
\[
\simeq \O(4p^n [O]),
\]
with the first two isomorphisms being canonical while the third we take to be given by
\[
(h_{b/p^n})^{p^n}.
\]

There is an isomorphism
\[
[2]^*\O(D_0)\simeq \O(D_0+D_1+D_2+D_3)\simeq \O(4D_0)
\]
with the last isomorphism being induced by the function $h_{b}\circ [p^n]$.

Finally, an isomorphism
\[
\O(4p^n [O])^{p^n}\simeq \O(4D_0)
\]
is induced by $f^{4}_{b/p^n}$. Given $x\in E_v\setminus E_v[2]$, choose $y$ such that
$p^ny=x$. Taking fibers above the points $y$ and $2y$, we obtain a commutative diagram 
\[
\xymatrix
@R=10ex @C=10ex
{
(\Hh_n)^4_y			\ar[r]^-{(\cdot)^{p^n}} 			\ar[d]^-{h^{-p^n}_{b/p^n}}_-\cong		&
(\Hh_n)^{4p^n}_y		\ar[r]^-{f^4_{b/p^n}}_-\cong	\ar[d]_-\cong^-{h^{-p^{2n}}_{b/p^n}}	&
\Oo(4D_0)_y											\ar[d]_-\cong^-{(h_b \circ p^n)\inv} 	\\
([2]^*\Hh_n)_y			\ar[r]^-{(\cdot)^{p^n}}			\ar[d]_-\cong							 &
([2]^*\Hh_n^{p^n})_y	\ar[r]^-{f_{2b/p^n}\circ [2]}_-\cong	\ar[d]_-\cong &
([2]^*\Oo(D_0))_y 										\ar[d]_-\cong	\\
(\Hh_n)_{2y}			\ar[r]^-{(\cdot)^{p^n}}	 &
(\Hh_n)_{2y}^{p^n} 	\ar[r]^-{f_{2b/p^n}}_-\cong &
\Oo(D_0)_{2y}
}
\]
where the lower vertical arrows are all natural base-change maps. The maps induced by
functions have all been based so as to make all diagrams commutative.
The maps are also clearly compatible with the action of the Galois group $G_v$.

\ssegment{d15}{}
We consider now the relation between the action of $g\in \cG_n$
and the composition of the leftmost vertical isomorphisms in the diagram, which we will denote by
\[
B(y):(\cH_n)_y^4\simeq (\cH_n)_{2y}.
\]
In the following, we will give the argument pointwise over $E$, even though
the underlying discussion is about the corresponding scheme isomorphism
\[
B:(\cH_n)^4\simeq [2]^*(\cH_n).
\]

\subsection*{Lemma}
There is a commutative diagram
\[
\xymatrix{
(\Hh_n)^4_y 				\ar[r]^-{g_1^4}_-\cong	 \ar[d]^-{B(y)}_-\cong &
(\Hh_n)^4_{\rho(g)(y)}		\ar[d]^-{B(\rho(g))(y)}_-\cong \\
(\Hh_n)_{2y}				\ar[r]^-{2g_1}_-\cong &
(\Hh_n)_{(2\rho(g))(2y)}
}
\]
where we denote by $2g_1$ the anti-symmetric lift of the element $2\rho(g)$.

\begin{proof}
We consider the isomorphism
\[
B(\rho(g)(y))\circ g_1^{\otimes 4} \circ B(y)^{-1}: (\cH_n)_{2y}\simeq (\cH_n)_{(2\rho(g))(2y)}
\]
lifting the action of $2\rho(g)$.
We need only check that
\[
i(B(\rho(g)(y))\circ g_1^{\otimes 4} \circ B(y)^{-1})=B(-y)\circ [(g_1)^{\otimes 4}]^{-1} B(-\rho(g)(y))^{-1} .
\]
For this, we embed the previous diagram into the bigger diagram
\[
\xymatrix{
(\Hh_n)^4_{-y}				\ar[r]^-{[-1]}_-\cong		\ar[d]^-{B(-y)}_-\cong	&
(\Hh_n)^4_y				\ar[r]^-{g_1^4}_-\cong 		\ar[d]^-{B(y)}_-\cong	&
(\Hh_n)^4_{\rho(g)(y)}		\ar[r]^-{[-1]}_-\cong		\ar[d]^-{B(\rho(g)(y))}_-\cong	&
(\Hh_n)^4_{-\rho(g)(y)}								\ar[d]^-{B(-\rho(g)(y))}_-\cong	\\
(\Hh_n)_{-2y}				\ar[r]^-{[-1]}_-\cong	&
(\Hh_2)_{2y}				\ar[r]^-{2g_1}_-\cong	&
(\Hh_n)_{(2\rho(g))(2y)}	\ar[r]^-{[-1]}_-\cong	&
(\Hh_n)_{-(2\rho(g))(2y)}.
}
\]
The two squares on the left and right are clearly commutative. But
$$[-1]\circ B(\rho(g)(y))\circ g_1^{\otimes 4} \circ B(y)^{-1}\circ [-1]$$
$$= B(-\rho(g)(y))\circ [-1]\circ g_1^{\otimes 4} \circ [-1]\circ B(-y)^{-1}$$
$$=B(-\rho(g)(y)) (g_1^{-1})^{\otimes 4} B(-y)^{-1}$$
$$=B(-\rho(g)(y)) (g_1^{\otimes 4})^{-1} B(-y)^{-1}$$
at every $y$ as desired.
\end{proof}

\ssegment{d16}{}
Choose $y$ as above so that $p^ny=x$ and let $v\in X_n$ lie above $y$.
Then for $\g\in G_v$, we have
\[
\ga(v) = c^x(\ga)v = c^x_2(\ga)c^x_1(\ga)v \in (\Hh_n)_{\rho(c^x(\ga))y}.
\]
Hence,
\[
\ga(v^{\otimes 4}) = (\ga(v))^{\otimes 4} = (c^x_2(\ga))^4(c_1^x(\ga)v)^{\otimes 4}.
\]
(Recall from segment \ref{3_a} that tensor powers restrict to ordinary powers in $\mu_{p^n}$.) Hence,
\[
\g(B(y)(v^{\otimes 4}))=B(\g(y))((\g v)^{\otimes 4})=B(\g(y))((c^x_2(\g))^4 (c^x_1(\g)v)^{\otimes 4})
\]
which by Lemma \ref{d15} 
\[
=(c^x_2(\g))^4(2c_1^x(\g))(B(y)(v^{\otimes 4})).
\]

\ssegment{d17}{}
We use the diagram of segment \ref{d14} to find that the map
\[
\cH_{2y}\rTo \O(D_0)_{2y}
\]
sends $B(y)(v^{\otimes 4})$ to
\[
(h_b(p^ny))^{-1} (1_{\O(D_0)})(2y)=(h_b(x))^{-1}(1_{\O(D_0)})(2y).
\]
Therefore,  an element in the inverse image $(X_n)_{2y}$ of $(1_{\O(D_0)})(2y)$
is
\[
(h_b(x))^{p^{-n}}B(y)(v^{\otimes 4}).
\]
Therefore, if we let $k(\cdot)_{p^n}$ denote the mod $p^n$ abelian Kummer map
\[
F_v^* \to H^1 \big( G_v,  \ZZ/p^n (1) \big)
\]
given in terms of the choice of a ($p^n$)\textit{th} root by
\[
k(a)_{p^n}(\ga) :=
  \frac{\ga \left( \sqrt[p^n]{a} \right)}{\sqrt[p^n]{a}},
\]
then the Galois action
on this element
is given by the cocycle
\[
k(h_b(x))_{p^n}(c^x_2)^4(2c^x_1).
\]
This must be the same as the action via $c^{2x} ,$ by the compatibility of
the big diagram with the action of $G_v$. As we take the limit over $n$, we get the equality
$$c^{2x}=k(h_b(x))(c^x_2)^4(2c^x_1)$$
at the level of $p$-adic cocycles.

\ssegment{d18}{}
One last modification is that this calculation has produced the
class   $$[\pip(2b,2x)]\in H^1(G_v, \pip(2b)),$$
which we need to shift back to $H^1(G_v, \pip)$ to get the class $\hat j(2x)$.
For this, we need to compose with the class of
$\pip(b,2b).$
We claim that this $\pip(b)$-torsor corresponds to the cohomology class
$k(2)$, where $k$ denotes the profinite abelian Kummer map
\[
F_v^* \to H^1(G_v, \Z_p(1)).
\]
To see this, let
\[
T_OX := T_OE_v \setminus \{0\}
\]
denote the punctured tangent space at the origin. There's an $F_v$-rational isomorphism of vector groups
\[
\AA^1 \xto{\sim} T_OE_v
\]
sending $1 \mapsto b$, hence an isomorphism of schemes
\[
\Gm \xto{\sim} T_OX
\]
which sends $1$ to $b$ and $2$ to $2b$. The theory of tangential fiber functors gives rise to an associated morphism of fundamental groupoids. In particular, there's a map
\[
\Zp(1) = \pi_1^{(p)}(\overline \GG_m,1) \to \pi_1^{(p)}(\overline X, b)
\surj \pip(b),
\]
and the induced map
\[
H^1(G_v, \Zp(1)) \to H^1(G_v, \pip(b))
\]
sends the torsor $\pi_1^{(p)}(1,2)$ to $\pip(b,2b)$. A straightforward calculation, carried out in \S14 of Deligne \cite{Deligne89}, shows that the former is represented by the Kummer cocycle $k(2)$ as claimed. 

Therefore,

\subsection*{Theorem}
Let $x \in X(F_v)$, let $c^x$ be an associated anticocycle
\[
G_v \to \varprojlim \Gg_n
\]
as in segment \ref{d13}, let 
\[
c^x = c_2^xc_1^x
\]
denote the decomposition of $c^x$ with $c^x_2$ taking values in $\Zp(1)$ and $c^x_1$ anti-symmetric (same segment), let $h_b$ denote the meromorphic function
\[
h_b = \frac{2y+a_1x + a_3}{2}
\]
of segment \ref{d14}, and let $k$ denote the Kummer map. Then
\[
k(2h_b(x))(c^x_2)^4(2c_1^x)
\]
is an anti-cocycle associated to the point $2x$.

\ssegment{d19}{}
We can now push out through the homomorphism $\pip\rTo U_2$.
\subsection*{Corollary}
Let
\[
j_v:\cX_v(F_v)\rTo H^1(G_v, U_2)
\]
be the unipotent Albanese map of level 2 at $v$. If
$j_v(x)=[c^x_2c^x_1]$, then 
\[
j_v(2x)=k(2h_b(x)) (c^x_2)^4(2c^x_1).
\]

\segment{phiv2x}{Proof of theorem \ref{Ner0}}

\ssegment{d22}{Lemma}
Suppose $a \in X(F_v)$ reduces to $O$ mod $m_v=(\pi_v)$ (the maximal ideal of $\O_{F_v}$). Then there exists an $a' \in E_v(F_v)$ such that
\[
a=p^na'.
\]

\begin{proof}
Let $D(O)$ denote the residue disk of $O$ inside $E_v(F_v)$. 
Referring to Silverman \cite{SilArEl}, Proposition 2.2 of Chapter VII, combined with Example 3.1.3 and Proposition 3.2 of Chapter IV, together provide a bijection
\[
D(O) \xto{\sim} m_v,
\]
plus a decreasing filtration of $D(O)$ by subgroups $D^i(O)$ compatible with the filtration on $m_v$ by powers, such that for each $i$, the induced map
\[
\tag{*}
D^i(O) / D^{i+1}(O) \xto{\sim} m^i_v/m^{i+1}_v
\]
is an isomorphism of groups. Moreover, $D(O)$ is separated and complete with respect to the filtration by the subgroups $D^i(O)$. 

Since $F_v$ contains $\QQ_l$, $l\neq p$, the group $m_v$ is $p$-divisible. We may use the group isomorphisms (*) to construct a Cauchy sequence $\{a_i\}$ in $D(O)$ with $p^n a_i \equiv a \mod D^i(O)$. Its limit $a'$ is a $p^n$\textit{th} root of $a$ as hoped.
\end{proof}

\ssegment{d22pp}{Lemma}
Suppose $a \in X(F_v)$ reduces to $O$ mod $m_v=(\pi_v)$. Then the anti-cocycle
\[
c^a: G_v \to \Gg_n(\overline F_v)
\]
associated to $a$ takes values in $\mu_{p^n}$.

\begin{proof}
According to lemma \ref{d22}, $a$ possesses an $F_v$-rational $p^n$\textit{th} root
\[
a' \in Y_n = E \setminus E[p^n]. 
\] 
Let $a''$ be an $\overline F_v$-point of $X_n$ lying above $a'$. For $\ga \in G_v$, the element
\[
g:=c^a(\ga) \in \Gg_n(\overline F_v)
\]
was defined in segment \ref{d13} by
\[
g(a'') = \ga(a''). 
\]
Recall from segment \ref{GnExt} that we have a commutative diagram like so:
\[
\begin{tikzcd}
&&
X_n \arrow[loop below, distance = 25]
\arrow{r}{s}
&
Y_n \arrow[loop below, distance = 25] 
\\
0 \arrow[r]
&
\mu_{p^n} \arrow[r]
&
\Gg_n \arrow{r}{\rho}
&
E_v[p^n] \arrow[r]
&
0.
\end{tikzcd}
\]
Since $\ga$ acts trivially on $a'$, we have 
\[
s(a'') = s(ga''),
\]
and because of the commutativity of the diagram, the latter equals
\[
\rho(g)(s(a'')).
\]
Since the action of $E[p^n]$ on $Y_n$ is free, it follows that $\rho(g) =0$, hence that $g \in \mu_{p^n}$.
\end{proof}

\ssegment{d23}{}
Near $O$, in the coordinate
$z$,  we have
\[
f_{b}=z^{1-p^{2n}}+z^{2-p^{2n}}g(z)
\]
with $g(z)\in \Oo_v[[z]]$. Also,
\[
z(a')=(1/p^n)z(a)+z(a)^2h(z(a)),
\]
for a power series $h\in \Oo_v[[z]]$. Therefore,
\[
f_{b/p^n}=(p^nz)^{1-p^{2n}}+(p^nz)^{2-p^{2n}}g_1(p^nz)
\]
with $g_1(t)\in \Oo_v[[t]]$ and
\[
f_{b/p^n}(a')=z(a)^{1-p^{2n}}+z(a)^{2-p^{2n}}H(z(a)),
\]
where $H(t)\in \Oo_v[[t]]$. Hence,
\[
f_{b/p^n}(a')=z(a)^{1-p^{2n}}u
\]
for a unit $u\equiv 1 \mod m_v$ and (since $c^a$ takes values in $\mu_{p^n}$, and since
\[
k(\cdot)_{p^n}: H^0(\Gm) \to H^1(\mu_{p^n})
\]
is a homomorphism to a $\ZZ/p^n$-module)
\[
c^a = k(f_{b/p^n}(a'))_{p^n}=k(z(a))_{p^n}.
\]
Taking the limit over $n$, we see that
in $H^1(G_v, \Z_p(1))$,
the class $j_v(a)$ is identified with the Kummer class of $z(a)$. 
Hence, for $a$ reducing to $O$ mod $m_v$, we have
\[
\phi_v (a)=\log \chi(rec_v(z(a)))= \log ( l^{v(z(a))})=-\log |z(a)|_v=(1/2)\log|x(a)|.
\]
By this formula, the function $\phi_v$ is bounded on the complement
\[
D(O) \setminus U
\]
of any open set $U$ containing $O$ inside the residue disk about $O$. On the other hand, by Kim-Tamagawa \cite[Corollary 0.2]{KimTamagawa},
\[
\phi_v(E_v \setminus D(O)) = \phi_v(X(\Oo_v))
\]
is finite. Thus, $\phi_v$ is bounded on the complement in $E_v$ of any $v$-adic neighborhood of $O$.
 
 \ssegment{d24}{}
 Finally, by corollary \ref{d19}, we have 
\begin{align*}
\phi_v(2a) &= (\log \chi) \cup j_v(2a)  \\
&= (\log \chi) \cup k(2h_b(a))(c_2^a)^4 	\\
&= (\log \chi) \cup k(2h_b(a)) + 4(\log \chi) \cup c_2^a \\
& =4\phi_v(a)-\log|2h_b(a)| \\
&=4\phi_v(a)-\log |(2y+a_1x+a_3)(a)|_v.
\end{align*}
This completes the proof of theorem \ref{Ner0}.


\segment{ran}{The range of a $p$-adic local N\'eron function}

\ssegment{Ner1}{}
We temporarily relax our assumption that $F_v$ is unramified over $\QQ_l$, and let $e$ denote the ramification degree. We normalize our absolute value $|\cdot|_v$ by $|l| = l\inv$. When taking $p$-adic logarithms of absolute values, we may artificially define
\[
\log (l^{n/e}) = \frac{n}{e} \log l.
\]
We also write $v = -\log |\cdot|$ (a valuation with values in the totally ordered subgroup $\ZZ \frac{\log l}{e}$ of $\Qp$), and we write
\[
\ord = \frac{e}{\log l}v.
\]

\subsection*{Proposition}
Suppose the function
\[
\la: E_v(F_v) \setminus \{O\} \to \Qp
\]
is a $p$-adic local N\'eron function in the sense of segment \ref{Ner0}.

(a) If $a \in E_v(F_v)$ reduces to a nonsingular point, then
\[
\la(a) = \max \{0, -\frac{1}{2}v(x(a)) \}.
\]

(b) Assume $E_v$ has multiplicative reduction and suppose $a\in E_v(F_v)$ reduces to a singular point. Let $N = \ord \Delta(E_v)$. Let $E_{v,0}(F_v)$ denote the group of points reducing to nonsingular points. We choose representatives $\{0, \dots, N-1\}$ for $\ZZ/N$. Then there is a unique isomorphism
\[
n: E_{v}(F_v) / E_{v,0}(F_v) = \ZZ/N
\tag{$*$}
\]
such that
\[
\la(a) = \frac{n(a)(N-n(a))}{2N^2}\log | \Delta |.
\tag{$**$}
\]

\ssegment{Ner2}{}
For the proof of proposition \ref{Ner1} we follow the treatment in chapter VI of Silverman \cite{SilAd}. In order to accord with the normalization used there, we set
\[
\la' = \la+\frac{1}{12}v(\Delta).
\]
Then $\la'$ satisfies (i), (ii), and
\begin{enumerate}
\item[(iii)'] For all $a \in E_v(F_v)$ with $[2]a \neq 0$, 
\[
\la'([2]a) = 4\la'(a) + v((2y+a_1x+a_3)(a)) - \frac{1}{4}v(\Delta).
\]
\end{enumerate} 
The proof of Theorem 4.1 of loc. cit. applies with the real logarithm replaced by the $p$-adic logarithm to show that
\[
\la'(a) = \frac{1}{2}\max\{v(x(a)\inv),0\} +\frac{1}{12}v(\Delta),
\]
which establishes (a). Our proof of (b) is similar; we nevertheless take the time to fill in some details in segments \ref{Ner3}--\ref{Ner4} below.

\ssegment{Ner3}{}
The proof of Theorem 1.1 of loc. cit. applies with the real logarithm replaced by the $p$-adic logarithm to show that properties (i)--(iii)' uniquely determine $\la'$. Let $L_v$ be a finite extension of $F_v$ of ramification degree $e'$ over $F_v$, suppose (b) has been established over $L_v$, and let $\la'$ be a function
\[
E_v(F_v) \setminus \{O\} \to \Qp
\]
 which satisfies properties (i)--(iii)'. Then the formulas given in parts (a) and (b) give us a function
 \[
 \la'_L: E_v(L_v) \setminus \{O\} \to \Qp
 \]
which, by uniqueness, extends $\la'$. Our preferred generator $a_0$ of $E_v(L_v)/ E_{v,0}(L_v)$ gives us a preferred generator $e'a_0$ of $E_v(F_v)/E_{v,0}(F_v)$. We have
\[
N = N_L / e'
\]
and for $a \in E_v(F_v)$ we set
\[
n(a) = n_L(a)/e'.
\]
Then
\[
\frac{n_L(N_L-n_L)}{2N_L^2}\log | \Delta_L | = \frac{n(N-n)}{2N^2}\log | \Delta |
\]
which establishes \ref{Ner1}($**$) over $F_v$. The uniqueness of \ref{Ner1}($*$) follows as in Lemma 5.1 of Silverman \cite{SilvermanComputing}. So after possibly replacing $F_v$ by a finite extension, we may assume $E_v$ has split multiplicative reduction.

\ssegment{Ner4}{}
It follows that $E_v$ is isomorphic to a Tate curve $E_q$ for some $q \in F_v^*$ with $|q| <1$ and $v(\Delta) = v(q)$. Let $\psi$ denote the induced map
\[
\psi: F_v^* \surj E_v(F_v).
\]
By Chapter V \S4 of Silverman \cite{SilAd}, 
\[
\psi : F_v^*/ q^\ZZ \xto{\cong} E_v(F_v)
\]
restricts to
\[
\Oo_v^* \xto{\cong} E_{0,v}(F_v).
\]
So the isomorphism
\[
E_v(F_v)/ E_0(F_v) \cong \ZZ/N
\]
is realized as
\[
\frac{F_v^*}{q^\ZZ \Oo_v^*} \xto{\ord} \ZZ/N.
\]
Thus, if
\[
a= \psi(u)
\]
with $0 < v(u) <v(q)$, we have
\[
\frac{n(N-n)}{2N^2}\log|\Delta| = \frac{1}{2}\left( \frac{v(u)^2}{v(q)} - v(u)   \right).
\]
So \ref{Ner1}($**$) is equivalent to 
\[
\la'(\phi(u)) = \frac{1}{2}B_2\left(\frac{v(u)}{v(q)} \right)v(q),
\]
where $B_2(T) = T^2-T +1/6$. We then set 
\[
\la'(\phi(u)) := \frac{1}{2}B_2\left( \frac{v(u)}{v(q)} \right)v(q) + v(\theta(u))
\]
where
\[
\theta(u) = (1-u)\prod_{m\ge 1} \frac{(1-q^m u)(1-q^m u\inv)}{(1-q^m)^2},
\]
and check that $\la'$ satisfies (i)--(iii)'. The proof of Chapter VI, Theorem 4.2 of Silverman \cite{SilAd} applies with the real logarithm replaced by the $p$-adic logarithm throughout.

This completes the proof of proposition \ref{Ner1}.

\segment{phirange}{Corollary}
Suppose $F_v$ is unramified over $\QQ_l$ and suppose $E_v$ has semistable reduction. Let $N_v=v(\Delta(E_v))$. Then the possible values for
$\phi_v$ on $X(\Oo_v)$ are
\[
-(n(N_v-n)/2N_v)\log l;\ \ \ 0\leq n< N_v.
\]

\begin{proof}
By proposition \ref{Ner1}, this follows from theorem \ref{Ner0}.
\end{proof}

\section{Punctured elliptic curves of low rank}
\label{S6}

\segment{d29}{}
We put ourselves in the situation and the notation ($\Ee$, $\Xx$, $\al$, $\be$, $b$, $S$, $N_l$, $W_l$, ...) of segment \ref{11_a} with $p$ an odd prime of good reduction, and $T = S \cup \{p\}$, and with the goal of proving the theorem stated there,
we begin by computing the image of
\[
\loc_p: \Sel^2(\Xx) \to H^1_f(G_p, U_2).
\]
We have the exact sequence
\[
0 \to \Q_p(1) \to U_2 \to V_p(E) \to 0,
\]
where $V_p(E)=T_p(E)\otimes \Q_p$ is the $\Q_p$-Tate module of $E=\cE\otimes \Q$. We recall that
\[
H^0(G_T, V_p(E)) = H^0(G_p, V_p(E)) = 0,
\]
so we have inclusions like so.
\[
\xymatrix{
H^1(G_T, \Qp(1)) \ar[r] \suphook[d]_-\theta &
H^1(G_p, \Qp(1)) \suphook[d] \\
H^1(G_T, U_2) \ar[r] &
H^1(G_p, U_2)
}
\]
It is straightforward to check that in this context, maps of Galois modules send crystalline classes to crystalline classes, so these inclusions induce inclusions like so.
\[
\xymatrix{
H^1_f(G_T, \Qp(1)) \ar[r] \suphook[d]_-{\theta_f} &
H^1_f(G_p, \Qp(1)) \suphook[d] \\
H^1_f(G_T, U_2) \ar[r] &
H^1_f(G_p, U_2)
}
\]

\segment{6_20}{Lemma}
If we assume, as in theorem \ref{11_a}, that $\cE(\Z)$ has rank zero\footnote{To avoid misunderstanding, we
remind the reader that $\cE$ refers to the compact curve, so that $\cE(\Z)=E(\Q)$, where
$E$ is the generic fiber of $\cE$. That is, what we write as $\cE(\Z)$ is what is usually called the rational points of $E$, while our $\cX(\Z)$ is sometimes confusingly referred to as  the integral points of $E$.} and that $\Sha_E[p^{\infty}] < \infty$, then $\Sel^2(\Xx)$ is contained in the image of $\theta_f$.

\begin{proof}
It is a general fact (which is straightforward to check) that the map 
\[
H^1(G_T, U_{n+1}) \to H^1(G_T, U_n)
\]
restricts to a map of Selmer schemes. On the other hand, 
we have an inclusion $\Sel^1(\Xx) \subset \Sel^1(\Ee)$, and
\[
\Sel^1(\Ee) = \Ee(\QQ)\otimes \Qp =0
\]
by corollary \ref{14a}, so
\[
\Sel^1(\Xx)=0.
\]
Hence each $P \in \Sel^2(\Xx)$ is $\theta(Q)$ for some $Q \in H^1(G_T, \Qp(1))$. To see that $Q$ is crystalline at $p$, we recall
 that
\[
H^0(G_p, V_p(E)\otimes B_\m{cr}) = 0,
\]
which implies that the map $\theta_B$ in the following diagram
\[
\xymatrix{
H^1(G_T, \Qp(1)) \ar[r] \suphook[d]_-\theta &
H^1(G_p, \Qp(1)) \suphook[d] \ar[r] &
H^1(G_p, B_\m{cr}(1) ) \suphook[d]^-{\theta_B} \\
H^1(G_T, U_2) \ar[r] &
H^1(G_p, U_2) \ar[r] &
H^1(G_p, U_2 \otimes B_\m{cr})
}
\]
is injective as shown, so that $\theta(Q)$ crystalline implies $Q$ crystalline.
\end{proof}

This allows us to regard $\Sel^2(\Xx)$ as a subset of $H^1_f(G_T, \Qp(1))$, and to compute its image in $H^1_f(G_p, \Qp(1))$.

\segment{11_b}{}
Recall (for instance from segment 6.2 of \cite{CKtwo}) that $H^1_f(G_T, \Q_p(1))$ can be realized as the subspace of $\Q^*\otimes_{\Z}\Q_p$ spanned by elements that are units outside $S$. Since $\ZZ^* \otimes \Qp =0$,
we have
\[
H^1_f(G_T, \Q_p(1))\simeq [\Z_S]^*\otimes_{\Z} \Q_p \inj  H^1_f(G_p, \Q_p(1)) \oplus \bigoplus_{v\in S} H^1(G_v, \Q_p(1)).
\]
We define the function 
$$\phi_v: H^1(G_v, \Q_p(1))\rTo \Q_p$$
by
$$c \mapsto \log\chi \cup c$$
(including $v=p$). We put these together to define
\[
\phi:  H^1_f(G_p, \Q_p(1))\oplus \bigoplus_{v\in S} H^1(G_v, \Q_p(1)) \to \Q_p
\]
by
\[
\phi( (c_v)_{v\in T})=\sum_v \phi_v(c_v).
\]

\segment{11_c}{Lemma}
For elements $c\in H^1_f(G_T,\Q_p(1))$, we have
\[
\phi \big( (\loc_v c )_{v\in T} \big)=0.
\]

\begin{proof}
We have 
\begin{align*}
\phi \big( (\loc_v c )_{v\in T} \big) 
	&= \sum_{v \in T} \log \chi \cup \loc_v(c) \\
	&= \sum_{v \in T} \loc_v (\log \chi \cup c) \\
	&= \sum_{\mbox{all places } v} \loc_v (\log \chi \cup c)
\end{align*}
since the contributions away from $T$ vanish.

By global class field theory (see, for instance, Tate \cite[Theorem B, \S11]{Tate}), the composite
\[
H^2(G_F, \Gm) \to \bigoplus_v H^2(G_v, \Gm) \to \QQ/\ZZ
\]
is equal to zero. By Hilbert's theorem 90, the cohomologies with $\mu_{p^n}$-coefficients inject into the cohomologies with $\Gm$-coefficients. Taking inverse limits and tensoring with $\Qp$, we find that the composite
\[
H^2\big(G_F, \Qp(1)\big)
\to \bigoplus_v H^2\big(G_v, \Qp(1) \big)
\to
\Qp
\]
is equal to zero, which completes the proof of the lemma.
\end{proof}

\segment{}{}
 For $l\neq p$,  we saw in corollary \ref{phirange} that
on $j_l\cX(\Z_l)$ the function $\phi_v$ takes  the values
$$-(n(N_l-n)/2N_l)\log l,$$
where $N_l=\ord_l \Delta_{\cE}$. As in the introduction, we define
\[
W_l := \set
{\frac{n(N_l-n)}{2N_l} \log l}
{0 \le n <N_l}
\]
and
\[
W := \prod_{l \in S} W_l,
\]
and for $w = (w_l) \in W$, we set
\[
|| w || := \sum_{l\in S} w_l.
\]
According to lemma \ref{11_c}, if $c \in \Sel^2(\Xx)$, we have
\[
\phi_p(\loc_p(c))=-\sum_{v\in S} \phi_v(\loc_v(c))= \| w\|
\]
for some vector $w \in W$.

\subsection*{Proposition}
With assumptions as above
we have
\[
\loc_p(H^1_{\Z}(U_2))=\bigcup_{w\in W} \set{ \eta \in H^1_f(G_p, U^3\bs U^2) }{ \phi_p(\eta)=\|w\| }.
\]

\begin{proof}
The inclusion $\subset$ has already been shown. To see that these equations define exactly the image,
note that the local reciprocity law
\[
\log \chi \cup k(a)=\log \chi (rec_v(a))
\]
for $a\in \Q_v^*$ shows that we get an isomorphism
\begin{align*}
\tag{*}
(\log \chi) \cup (\,\cdot\,) :
H^1_f(G_p, \Q_p(1)) \oplus
	& \bigoplus_{v\in S} H^1(G_v, \Q_p(1))
\\
	& \simeq H^2(G_p, \Q_p(1))\oplus 
	\bigoplus_{v\in S} H^2(G_v, \Q_p(1))
\\
	&\simeq \bigoplus_{v\in T} \Q_p.
\end{align*}

Indeed, for $v \neq p$,
$
H^1(G_v, \Qp(1))
$
is one dimensional and generated by the class of $k(v)$. We see this by the exact sequence
\[
0 \to \ZZ_v^* \to \QQ_v^*\to \ZZ \to 0
\]
and the fact that the kernel has to map to zero under the Kummer map (since $v \neq p$). So it suffices to show that
\[
\log(\chi(\rec(v)))
\]
is non-zero in $\Qp$. But
\[
\chi(\rec(v))=\chi(\Fr_v)
\]  
is just $v \in \ZZ_p^*$, an element of infinite order. So its log is non-zero.

For $v=p$, $H^1(G_p, \Qp(1))$ is two-dimensional. But we've already discussed the fact that $H^1_f(G_p, \Q_p(1))$ is one-dimensional, generated by the Kummer image of the units in $\Zp$. Thus, it suffices to show that $\chi(\rec(\ZZ_p^*))$ is of infinite order. But in fact, $\chi(\rec(\cdot ))$ just induces an isomorphism $\ZZ_p^*\simeq \Aut(\Zp(1))$ by the definition of the reciprocity map in local class field theory (\cite[p. 146]{Serre}).

The isomorphism (*) will take the (injective) image of
$H^1_f(G_T, \Q_p(1))$ to the (injective) image of $H^2(G_T, \Q_p(1))$, which is exactly equal to the kernel of the sum map
$$\bigoplus_{v\in T} \Q_p\xto{\Sigma_v}{\Q_p}.$$
Since
$$\dim H^1_f(G_T, \Q_p(1))=\dim [(\Z_S)^*\otimes \Q_p]= |T|-1,$$
we see thereby that $(\log\chi) \cup (\,\cdot \,)$ takes $H^1_f(G_T, \Q_p(1))$ also isomorphically to the kernel of the sum map.
On the other hand, we have seen that
$$\log\chi \cup: \bigoplus_{v\in S} j_v(\cX(\Z_p))\simeq \bigoplus_{v\in S} W_v\subset \bigoplus_{v\in S} \Q_p.$$
Therefore, the subspace
$$H^1_{\Z}(G, U_2)\subset H^1_f(G_T, \Q_p(1)),$$
which is defined as the inverse image of $\bigoplus_{v\in S} j_v(\cX(\Z_p))$, is exactly defined by 
\[
\bigcup_w \set { \eta\in H^1_f(G_T, \Q_p(1)) }{ \sum \phi_p(\loc_p(\eta))=\|w\|}.
\]
Hence, the $p$-component of elements of $H^1_{\Z}(G, U_2)$, that is, its image under $\loc_p$, is
exactly defined by the equations in the statement of the proposition.
\end{proof}

\segment{11_d}{Remark}
According to proposition \ref{11_c}, the equality
\[
\cX(\Z)=\cX(\Z_p)_2
\]
may be viewed as an exactness statement for the sequence 
\[
1 \to \cX(\Z) \to \prod_{v\in T} \cX(\Z_v) \xto {h_p} {\Q_p},
\]
in a manner reminiscent of class field theory. Since the map $h_p$ is quadratic, exactness here should be understood in the sense of pointed sets. Of course, this cannot hold literally, since we could take the image of an integral point in
$\prod_{v\in T} \cX(\Z_v)$ and move it inside its residue disk just at one $v\neq p$ without changing the height. The formula
\[
\cX(\Z)=\cup_w \Psi(w)\subset \cX(\Z_p),
\]
`the projection to the $p$-component', is one  recasting of this exactness that absorbs the ambiguity.

\segment{d30}{}
We can now prove theorem \ref{11_a}.
We follow the notation of \cite{KimMassey}, section 3 and let
\[
\Exp:L^{DR}_2/F^0\simeq H^1_f(G_p, U_2)
\]
denote the non-abelian Bloch-Kato exponential map from the Lie algebra of
the de Rham fundamental group. Recall that $L_n^{DR}$ may be realized as the quotient of the tensor algebra $T^\cdot H_1^{DR}(\Xx_\Qp)$
 modulo the $(n+1)^\m{st}$ power of the augmentation ideal. We denote by $A$, $B$ the elements of $L_2^{DR}$ associated to the basis of $H_1^{DR}(\Xx_\Qp)$ dual to $\{\al, \be\}$. According to lemma 3.2 of \cite{KimMassey}, $L_2^{DR}/F^0$ has basis
 \[
 \{A, A^2, AB, BA\}.
 \] 
According to the proof of corollary 0.2$'$ of \cite{BalakAppendix}, the map
$$j_p:\cX(\Z_p) \to H^1_f(G_p, U_1)$$
is given by 
$$j_p(z)=\log(z)\Exp(A)=(\int_b^z \a)\Exp(A),$$
while
$$j_p:\cX(\Z_p) \to H^1_f(G_p, U_2)$$
is given by
$$j_p(z)=(\log(z)\Exp(A), D_2(z)\Exp([A,B])),$$
where 
$$D_2(z)=\int_b^z \a \b.$$
By
 Proposition 3.3 of \cite{KimMassey}, we have
$$\phi_p(\Exp([A,B])=1.$$
Therefore,
$$\phi_p(j_p(z))=D_2(z).$$

From this, we see that 
$$\cX(\Z_p)_1=\cE(\Z_p)(tor)\setminus O.$$
For small primes $p$, it will  happen frequently that
$$\cX(\Z)=\cX(\Z_p)_1,$$
since global torsion on $\cE$ will often be equal to the local torsion for $p$ small.
But of course, this fails for large $p$, and one must look at level 2, which then imposes 
on $\cX(\Z_p)_2$ the pair of
conditions
 $$ \log(z)=0, \ D_2(z)=\|w\| $$
for some $w$, as in the statement of the theorem.

\segment{d31}{}
So far, we have tested conjecture \ref{13i} using the prime $p=5$
for 256 semi-stable elliptic curves of rank zero from Cremona's table, and found
$$\cX(\Z)=\cX(\Z_p)_2$$
for each of them. To give a rough sense of the data computed using the methods of \cite{BalakIterated}, the details of which can be found on \cite{sage2}, we present here a small table illustrating some
of the large $\|w\|$-values that come up as we go through the list.

\medskip
\begin{minipage}{0.5\textwidth}
    \begin{tabular}{| c | c |}
    \hline
Cremona & number of \\ label &  $||w||$-values \\
\hline
1122m1 &128\\
1122m2 & 384\\
1122m4 &84\\
1254a2 &140\\
1302d1 &72\\
1302d2 &96\\
1426b4 &64\\
1506a2 &112\\
1806h1 &120\\
2397b1 &72\\
2418b2 &64\\
2442h1 &78\\
2442h2 &84\\
2478c2 &68\\
2706d2 &120\\
2967c1 &72\\
2982j1 &160\\
2982j2 &140\\
3054b1 &108\\
\hline
          \end{tabular}
  \end{minipage}\qquad
  \begin{minipage}{0.5\textwidth}
    \begin{tabular}{| c | c |}
    \hline
Cremona  & number of \\ label &  $||w||$-values \\
\hline
3094d1 &72\\
3486o1&72\\
3774f1 &120\\
4026g1 &90\\
4134b1 &90\\
4182h1 &300\\
4182h2 &64\\
4218b1 &96\\
4278j1 &90\\
4278j2 &100\\
4434c1 &210\\
4514d1 &64\\
4602b1 &64\\
4658d2 &66\\
4774e1 &224\\
4774e2 &192\\
4774e3 &264\\
4774e4 &308\\
4862d1 & 216\\
\hline
       \end{tabular}
  \end{minipage}
\medskip

Hence, for example, for the curve $1122m2$, 
$$ y^2 + xy = x^3 - 41608x - 90515392$$
there are 384 of the $\Psi(w)$'s that potentially make up $\cX(\Z_p)_2$. Of these,
all but 4 end up being empty, while the points in those $\Psi(w)$ consist exactly of the integral points 
$$
(752 , -17800 ), (752, 17048 ), (2864 , -154024 ), (2864, 151160 ).
$$

\segment{d32}{}
Another kind of test  is to fix a few curves and let $p$ grow. For example, for the curve (`378b3')
$$y^2+xy=x^3-x^2-1062x+13590,$$
we found that
$$\cX(\Z_p)_2=\cX(\Z)=\{ (19, -9), (19, -10)\}$$
for $5\leq p\leq 97.$
As one might expect, as $p$ gets large, $\cX(\Z_p)_1$ becomes significantly larger than $\cX(\Z)$. For $p=97$,
we have
$$|\cX(\Z_{97})_1|=89.$$
However, imposing the additional constraint exactly cuts out the integral points for each $p$.

\segment{d33}{}
Is it conceivable that 
$$\cX(\Z)_{tor}=\bigcup_w \Psi_w\subset \cX(\Z_p)$$
even when $\cE$ has higher rank?  There are rather obvious relations with the conjecture on non-degeneracy of the $p$-adic height \cite{MazurSteinTate},
which we hope to investigate in a later work. At the moment, we have a small bit of evidence, having tested the equality numerically for
10 curves of rank one, 2 curves of rank two, and one curve of rank three. Some of the cases are quite dramatic, such as
Cremona label `82110bt2', which has rank one. In this case, there are 2700 different $\|w\|$-values to consider, each contributing
some  $\Psi(w)$. However, the only non-empty ones are those that contain the 14 integral torsion points
\cite{sage2}. 

\segment{d34}{}
We close this section with a brief mention of the framework for conjecture \ref{13i} when $\cE(\Z)$ has rank one, leaving a systematic
treatment to a later paper. As above, $\Ee$ denotes the regular minimal model of an elliptic curve over $\QQ$. 
Assume that there is a point
$y\in \cX(\Z)$ of infinite order.
In the case where the Tamagawa number of $\cE$ is 1, we  saw in \cite{KimMassey, BalakAppendix} that
\[
\loc_p:H^1_{\Z}(U_2) \to H^1_p(G_p, U_2)=\A^2
\]
is computed to be
\[
\A^1\to \A^2;
\]
\[
t\mapsto (t, ct^2),
\]
where
\[
c=D_2(y)/\log^2(y).
\]
So the image is defined by $x_2-cx_1^2=0$.
Meanwhile,
\[
\cX(\Z_p) \to H^1_f(G_p, U_2)
\]
is
\[
z\mapsto ( \log(z), D_2(z)).
\]
Thus,
\[
\cX(\Z_p)_2
\]
is the zero set of
\[
D_2(z)-c\log^2(z).
\]
It is sometimes convenient to write this defining equation as
\[
\frac{D_2(z)}{\log^2(z)}=c.
\]
In the earlier paper \cite{BalakAppendix}, we checked in a number of cases that the integral points
do indeed fall into the zero set. However, it turns out that 
\[
\cX(\Z)\subsetneq \cX(\Z_p)_2
\]
in the majority of cases, underscoring the importance of going up another level. (In fact, a superficial guess
based on the rank zero case would indicate that if the localization
\[
H^1_{\Z}(G, U_n) \rTo H^1_f(G_p, U_n)
\]
becomes injective at level $n$, then the conjecture might hold at level $n+2$.)

\segment{d35}{}
When we assume that $\mathcal{E}$ is semistable but with arbitrary Tamagawa numbers as above, the precise form of the equation becomes a bit delicate, and we will leave  a systematic treatment to a later paper. However, if we define in this case $$c = \frac{h(y)}{(\log(y))^2},$$
where $h$ is the $p$-adic height \cite{MazurSteinTate} (except that our convention for the height function multiplies theirs by $p$) and $y$ is a point of infinite order, then we can prove the following:

\subsection*{Proposition}
$$\mathcal{X}(\mathbb{Z})\subset  \bigcup_{w \in W} \{z \in \mathcal{X}(\mathbb{Z}_p) \large\mid D_2(z) + \frac{C}{2}(\log(z))^2 - c (\log(z))^2 =  \lvert\lvert w \rvert\rvert \},$$ where
\[
C = \frac{a_1^2 + 4a_2}{12} -\frac{\mathbf{E}_2(E,\alpha)}{12}
\]
is the special value of the $p$-adic modular form $\mathbf{E}_2$ associated to the pair $(E, \alpha)$.

  The equations on the right hand side should in fact define $\mathcal{X}(\mathbb{Z}_p)_2$, but we will not check this at the moment.

\begin{proof}Let $h_p(z) := h_p(z-O,z-O)$ denote the local height at $p$ of $z$. We first show that
\[
h_p(z) = D_2(z) + \frac{C}{2}(\log(z))^2.
\]
To do this, we use an interpretation of the local height at $p$ in terms of Coleman integrals as in Theorem 4.1 of \cite{BalBesColGr}. Note the following normalization: our global height (and local heights) are precisely half those in \cite{BalBesColGr}. 

In terms of our normalization of local heights, Theorem 4.1 of \cite{BalBesColGr} gives that $$h_p(z-O,z-O) = -\int_{b}^z \alpha \eta, $$ with  $[\eta] \cup [\alpha] = 1$. Note that $\eta$, which is found in the course of proving Corollary 4.2 of \cite{BalBesColGr}, is given by
\[
-\eta  = \eta_0 + C\alpha,
\]
where $\eta_0 =\beta$ and
\[
C = \frac{a_1^2 + 4a_2}{12} -\frac{\mathbf{E}_2(E,\alpha)}{12}
\]
is the special value of the $p$-adic modular form $\mathbf{E}_2$ associated to the pair $(E, \alpha)$.
 Then substituting appropriately, we have
\begin{align*} h_p(z) &= -\int_{b}^z \alpha \eta \\
&= \int_b^z \alpha \eta_0 + C\int_b^z \alpha \alpha\\
&= D_2(z) + \frac{C}{2}(\log(z))^2.\end{align*}

Thus we have that
\[
h_p(z) = D_2(z) + \frac{C}{2}(\log(z))^2.
\]
Finally, since
\[
c =\frac{h(y)}{(\log(y))^2} = \frac{h(z)}{(\log(z))^2},
\] 
we have 
\begin{align*}c(\log(z))^2 &= h(z) \\
&= D_2(z) + \frac{C}{2}(\log(z))^2 + \sum_{v \neq p} h_v(z-O, z-O),\end{align*}
and noting that the possible values of $-h_v(z-O, z-O)$ on integral points $z$ are precisely given by $\phi_v$, we conclude that 

$$\mathcal{X}(\mathbb{Z})\subset \bigcup_{w \in W} \{z \in \mathcal{X}(\mathbb{Z}_p) \large\mid D_2(z) + \frac{C}{2}(\log(z))^2 - c (\log(z))^2 =  \lvert\lvert w \rvert\rvert \}.$$\qedhere \end{proof}

\section{The thrice punctured line}
\label{S5}

\segment{d25}{}
Let $\cX=\P^1\setminus \{0,1,\infty\}$ and take $b$ to be the standard tangential base-point $\overrightarrow{01}$ based at $0\in \P^1$. For the basic facts here, we refer to \cite{kimi}.
We have
\[
U_1=\Q_p(1)\times \Q_p(1)
\]
and
\[
U^2/U^3=\Q_p(2),
\]
so there is
an exact sequence
\[
0 \to \Q_p(2) \to U_2 \to \Q_p(1)\times \Q_p(1) \to 0.
\]
The diagram
\[
\xymatrix{
\cX(\Z) \ar[d] \ar[r] & \cX(\Z_p) \ar[d] \\
H^1_{\Z}(G, U_1)  \ar[r]^-{\loc_p} & H^1_f(G_p, U_1)
}
\]
thus
becomes
\[
\xymatrix{
\emptyset \ar[d] \ar[r] & \cX(\Z_p) \ar[d] \\
0  \ar[r]^-{\loc_p} & \A^2.
}
\]
The map $j_p$ takes the form
\[
z\mapsto (\log(z), \log(1-z)),
\]
(see, for instance, Proposition 7.3 of Dan-Cohen--Wewers \cite{CKtwo}) so that
$\cX(\Z_p)_1$ is the common zero set of $\log(z)$ and $\log(1-z)$.
Since $z$ and $1-z$ must both be roots of unity, the only common zero possible is
$z=\z_6$ for a primitive sixth root of unity.
If $p=3$ or $p\equiv 2 \mod 3$, then $\z_6\notin \Q_p$, so that
$\cX(\Z_p)_1=\phi$. That is,
\subsection*{Proposition}
Conjecture \ref{13i} is true for $n=1$ when $\cX=\P^1\setminus \{0,1,\infty\}$
and $p=3$ or $p\equiv 2 \mod 3$.

\segment{d26}{}
When $p\equiv 1 \mod 3$
\[
\cX(\Z)=\phi \subsetneq \{\z_6, \z_6^5\} =\cX(\Z_p)_1
\]
and we must go to a higher level.
We have
$$H^1_f(G_p, U_2)=\A^3$$
and
$$j_p:\cX(\Z_p)\rTo \A^3$$
is given by
$$j_p(z)=( \log(z), \log(1-z), -Li_2(z)),$$
where $$Li_2(z)=\sum_n \frac{z^n}{n^2}=\int_b^z (dt/t) (dt/(1-t))$$
 is the {\em $p$-adic dilogarithm} (loc. cit.). 
Meanwhile, we still have
$$\Sel^1(\Xx)=0,$$
since $H^1(G_T, \Q_p(2))=0$ by Soul\'e's vanishing theorem \cite{soule}. Therefore, 
$$\cX(\Z_p)_2=\{z\ | \ \log(z)=0, \log(1-z)=0, Li_2(z)=0\},$$
and the question of whether $\cX(\Z)=\cX(\Z_p)_2$ reduces to checking
if $Li_2(\z_6)$ can be zero.  Note that (\cite{coleman},  Prop. 6.4)
$$Li_2(\z_6)+Li_2(\z_6^{-1})=-\log^2(\z_6)/2=0,$$
so that we need only discuss non-vanishing at one of the sixth roots.
This question was raised by Coleman in \cite{coleman}, page 207, remark 3.

\segment{d27}{}
We have checked numerically thus far that
$$Li_2(\z_6)\neq 0$$
for $p$ in the range $3\leq p <10^5.$
This may be carried out as follows. We define power series $g_n \in \Q[[v]]$ recursively by
\begin{align*}
&g_0(v) = v-1 - \frac{(v-1)^p}{v^p-(v-1)^p} \\
\intertext{and for $n \ge 1$,}
&g'_{n+1}(v) = -v^{-1}g_n(v)(1+v+v^2+\cdots) \\
&g_{n+1}(0)=0 .
\end{align*}
Then 
\[
Li_2(\zeta) = \frac{p^2}{p^2-1} g_2 ( \frac{1}{1-\zeta} )
\]
for any $(p-1)^{st}$ root of unity $\zeta$. Indeed, this is a special case of Propositions 4.2 and 4.3 of \cite{Lip}. Hence, since $ 1/(1 - \zeta_6) = \zeta_6$, it suffices to check that $g_2(\zeta_6) \neq 0$.
In fact, (for $p$ as above) $g_2(\zeta_6)$ is nonzero modulo $p$. Moreover, $g_2$ reduces modulo $p$ to a polynomial of degree $p-2$ which is determined by the reductions modulo $p$ of the same equations as above. So the verification may be performed rapidly with any computational software. We used Sage \cite{sage1}.

\segment{d28}{}
However, this may fall short of providing definitive evidence that
\[
Li_2(\zeta_6) \neq 0
\]
for all primes congruent to $1$ mod $3$: if for each $p$, the value of $g_2(\zeta_6)$ modulo $p$ were merely \textit{random}, the probability of $g_2(\zeta_6)$ being nonzero for $p$ in our range would be about 0.413. On the other hand, the probability that $g(\zeta_6)=0$ for some $p\equiv 1 \mod 3 $ would be 1. The point is that the product
$$\prod_{p=1 \mod 3} (1-1/p)$$
tends to zero, representing the probability that $g_2(\zeta_6)$ does  
not vanish mod $p$ for all $p\equiv 1 \mod 3$ (assuming these values are random,  
independently and evenly distributed variables). 

Suppose we want to falsify the randomness hypothesis. To do this, we  
could show that $g_2(\zeta_6) $ doesn't vanish for all $p< N$ for some  
large value of $N$. 
Unfortunately, the convergence of the product is extremely slow:
for $N=100,000 $ it is just 0.413, which 
does not give convincing evidence. This computation took several hours.  
Doing the computation up to $p< 10^6$ would take several days, and  
the probability would be 0.3775, which is not so much better.

\section{Remarks on curves of higher genus}
\label{S7}

\segment{12_a}{}
Let $\Xx \to \Spec \ZZ$ be the regular minimal model of a proper smooth curve of genus $\ge 2$ (case 2 of the trichotomy of segment \ref{13a}). We fix a base point $b \in \Xx(\ZZ)$. Then the associated map
\[
j_p:\cX(\Z_p)\rTo H^1_f(G_p, U_1)
\]
can be identified with the map
\[
\cX(\Z_p)\subset X(\Q_p) \inj J_X(\Q_p)\rTo T_eJ_X,
\]
where $J_X$ is the Jacobian of $X$ and $T_eJ_X$ is its tangent space
at the origin \cite{BlochKato}. 

\segment{12_b}{}
If we assume that
$\Sha_{J_X}[p^{\infty}]<\infty$, then it follows that
\[
J_X(\ZZ) \otimes \Qp = \Sel^1(\Xx).
\]
Indeed, the corollary of segment \ref{14a} applies equally to the abelian variety $J_X$:  
\[
J_X(\ZZ) \otimes \Qp = \Sel^1(J_X).
\]
Since the map $U_1(\Xx) \to U_1(J_X)$ is an isomorphism, we have a natural inclusion
\[
\Sel^1(\Xx) \subset \Sel^1(J_X).
\]
Conversely, if
\[
P\in \Sel^1(J_X)
\]
is an arbitrary Selmer class for the Jacobian,
then
\[
\loc_vP = 0 \mbox{ for all } v \neq p
\]
whence
\[
P \in \Sel^1(\Xx).
\]

\segment{6_13_a}{}
If we assume moreover $J_X(\ZZ)$ has rank zero, then it follows from segment \ref{12_b} that $\Sel^1(\Xx) =0$. So by segment \ref{12_a} we have
\[
\Xx(\Zp)_1 = \Xx(\Z_p) \cap J_X(\Zp)(tor).
\]

\segment{d36}{}
We apply this to the Fermat curve
\[
X_l: x^l+y^l=z^l
\]
for prime $l\geq 5$. It is a theorem of Coleman, Tamagawa, and Tzermias \cite[Theorem 2]{Tzermias} that
$$X_l(\Z_p)\cap J_{X_l}(\Z_p)(tor)$$
must satisfy $xyz=0$. Therefore, if $\z_l\notin \Q_p$, then
$$X_l(\Z)=X_l(\Z_p)\cap J_{X_l}(\Z_p)(tor).$$
(Using also the theorem of Wiles.) Meanwhile,
we have
\[
\rank J_{X_l}(\Z)=0
\]
for $l=5, 7$ \cite{Fadeev}. So by segment \ref{6_13_a}, we have

\subsection*{Proposition}
With notation as above, assume
\[
\Sha_{J_{X_l}}[p^{\infty}]<\infty.
\]
If $\cX_l$ is the minimal regular model of $X_l$, we have
\[
\cX_l(\Z)=\cX_l(\Z_p)_1
\]
for $l=5,7$ and $p\not\equiv 1 \mod l$. That is, conjecture \ref{13i} is true at level 1 for these $l$ and $p$.

\segment{d37}{}
It should be interesting to investigate conjecture \ref{13i} in relation to the many known results about torsion packets on curves, for example, for Fermat curves or modular curves \cite{BakerRibet}. It seems reasonable to  suspect that there should be more instances where
$\cX(\Z)=\cX(\Z_p)_1$ even when the Jacobian has positive rank. Another way to say this is that the classical method of Chabauty is usually applied with one choice of a differential form. The question here raised by conjecture \ref{13i} is how often the common zero set of {\em all} available abelian integrals will give us exactly $\cX(\Z)$.

\segment{d38}{}
In fact, it is relatively easy to produce example of affine curves of higher genus and good reduction primes $p$ for which the conjecture holds. We illustrate this by way of an example. Consider the elliptic curve $E$ with affine model
$$y^2=x^3-891x+4374.$$
It turns out that
$$E(\Q)\simeq \Z/4$$
and that the two points $P$ and $Q$ of order four are
$(-9, \pm 108).$
Thus, by making the substitution
$x+9=dt^2$, we get a cover
$$f:X\to E\setminus \{O\}$$
ramified exactly over $P$ and $Q$. This affine hyperelliptic curve has equation
$$y^2=d^3x^6-27d^2x^4-648dx^2+11664.$$
For any $p$ of good reduction  and $d$ such that $d$ is not a square in $\Q_p$, the point of order two
$(27,0)$ will not lift to $\cX(\Z_p)$. On the other hand, we have the commutative diagrams
\[
\xymatrix{
\cX(\Z) \ar[r] \ar[d]_-f^-\cong & \cX( \Z_p) \ar[d] \\
[\cE\setminus O ] (\Z) \ar[r] &  [\cE\setminus O](\Z_p)
}
\]
\[
\xymatrix{
\cX(\Z_p) \ar[d] \ar[r] &  H^1_f(G_p, \pi_1^{\Q_p}(\bX, b)) \ar[d] \\
[\cE\setminus O ](\Z_p) \ar[r] &H^1_f(G_p, \pi_1^{\Q_p}(\bE\setminus O , b))
}
\]
and
\[
\xymatrix{
H^1_{\Z}(G, \pi_1^{\Q_p}(\bX, b)) \ar[r] \ar[d] & 
 H^1_f(G_p, \pi_1^{\Q_p}(\bX, b)) \ar[d]  \\
H^1_{\Z}(G, \pi_1^{\Q_p}(\bE\setminus O, b)) \ar[r] &
 H^1_f(G_p, \pi_1^{\Q_p}(\bE\setminus O, b)).
}
\]
These imply that
$$\cX(\Z_p)_2\subset f^{-1}[\cE\setminus O](\Z_p)_2.$$
Hence, for such $d$ and $p$, it is easy to deduce that whenever
$$[\cE\setminus O] (\Z)=[\cE\setminus O](\Z_p)_2,$$
we also get
$$\cX(\Z)=\cX(\Z_p)_2$$
for free.

We can check this, for example, for $d=2$ and all
$3\leq p\leq 53$ such that $p\equiv 3 \mod 8$ and $p\equiv 5 \mod 8$.

\section{Remarks on $S$-integral points}
\label{S8}

\segment{d39}{}
Let $\Xx \to \Spec \ZZ$ denote a regular $\ZZ$-model of a hyperbolic curve over $\QQ$, and let $b$ denote a \emph{base point} as in segment \ref{13a}. As above we denote by $U_n$ the level-$n$ quotient of the unipotent $p$-adic \'etale fundamental group of $\Xx_{\bar\QQ}$ at $b$. We also use the same notation for the fundamental group of $\Xx_{\bar \QQ_p}$.  Let $S$ denote a finite set of primes of $\ZZ$, $p$ a prime of good reduction not in $S$, and let $T$ be a finite set of primes containing $S$ and $p$, as well as all primes of bad reduction. We define the \emph{$S$-integral Selmer scheme of $\Xx$ in level $n$} by
\[
\Sel^n_S(\Xx) = \bigcap _{v\neq p, v\notin S} \loc_v^{-1}[Im(j_v)]\subset H^1_{f}(G_T, U_n).
\]
This gives rise to a filtration
\[
\cX(\Z_p)_{S,n}:=j_p^{-1}(\loc_p(H^1_{\Z,S}(U_n)))
\]
on $\cX(\Z_p)$, which one might conjecture to converge to $\cX(\Z_S)$.

\segment{d40}{}
As an (admittedly small) step in this direction, we consider the case 
$\cX=\P^1\setminus \{0,1,\infty\}$ and $S=\{2\}$. Then the diagram
\[\xymatrix{
\cX(\Z[1/2]) \suphook[r] \ar[d]_-j &
\cX(\Z_p) \ar[d]^-{j_p} \\
\Sel^2_S(\Xx) \ar[r]_-{\loc_p} &
H^1_f(G_p, U_2) }
\]
becomes
\[\xymatrix{
\{2, 1/2, -1\} \suphook[r] \ar[d]_-j  & 
\cX(\Z_p) \ar[d]^-{j_p} \\
\A^2 \ar[r]_-{\loc_p} & 
\A^3 }
\]
where \cite{CKtwo}
\[
\loc_p(x,y)=((\log2) x, (\log2 )y, (1/2)(\log^22 ) xy ).
\]
Recall that
\[
j_p(z)=(\log(z), \log (1-z), -Li_2(z)),
\]
so that
$\cX(\Z_p)_2$
is the zero set of
\[
2Li_2(z)+\log (z) \log(1-z).
\]

The fact that $\{2, -1, 1/2\}$ is in the zero set, which we deduce here from the commutativity of the localization diagram for Selmer schemes, was noticed earlier by Coleman to be
 a consequence of standard
dilogarithm identities  (\cite{coleman}, remark on page 198). 
We have checked numerically that this is {\em exactly} the zero set
for $p=3,5,7$, so that
$$\cX(\Z[1/2])=\cX(\Z_p)_{\{2\},2}$$
in that case. The equality starts failing for larger $p$. Coleman already noted this failure for
$p=11$, since $\frac{-1\pm \sqrt{5}}{2}$,  for example, is in the zero set.
This fact, as well as the considerations of the previous sections, indicate the importance of investigating systematically
weakly global points of higher level.

\bibliography{eferences}

\begin{thebibliography}{MvdGE}

\bibitem[Bal1]{sage2}
Jennifer Balakrishnan.
\newblock Data page.
\newblock \url{http://math.harvard.edu/~jen/data.html}.

\bibitem[Bal2]{BalakIterated}
Jennifer~S. Balakrishnan.
\newblock Iterated coleman integration for hyperelliptic curves.
\newblock In {\em ANTS-X: Proceedings of the Tenth Algorithmic Number Theory
  Symposium, Open Book Series 1}. Mathematical Sciences Publishers, 2013.

\bibitem[BB]{BalBesColGr}
Jennifer~S. Balakrishnan and Amnon Besser.
\newblock {C}oleman--{G}ross height pairings and the p-adic sigma function.
\newblock {\em Journal für die reine und angewandte Mathematik (Crelle's
  journal)}.
\newblock To appear.

\bibitem[BBK]{BalakBK}
Jennifer~S. Balakrishnan, Robert~W. Bradshaw, and Kiran~S. Kedlaya.
\newblock Explicit {C}oleman integration for hyperelliptic curves.
\newblock In {\em Algorithmic number theory}, volume 6197 of {\em Lecture Notes
  in Comput. Sci.}, pages 16--31. Springer, Berlin, 2010.

\bibitem[BBM]{BalBesMul}
Jennifer~S. Balakrishnan, Amnon Besser, and Steffen M\"uller.
\newblock p-adic height pairings and integral points on hyperelliptic curves.
\newblock {\em Journal für die reine und angewandte Mathematik (Crelle's
  journal)}.
\newblock To appear.

\bibitem[BdJ]{Lip}
Amnon Besser and Rob de~Jeu.
\newblock {${\rm Li}^{(p)}$}-service? {A}n algorithm for computing {$p$}-adic
  polylogarithms.
\newblock {\em Math. Comp.}, 77(262):1105--1134, 2008.

\bibitem[BK]{BlochKato}
Spencer Bloch and Kazuya Kato.
\newblock {$L$}-functions and {T}amagawa numbers of motives.
\newblock In {\em The {G}rothendieck {F}estschrift, {V}ol.\ {I}}, volume~86 of
  {\em Progr. Math.}, pages 333--400. Birkh\"auser Boston, Boston, MA, 1990.

\bibitem[BKK]{BalakAppendix}
Jennifer~S. Balakrishnan, Kiran~S. Kedlaya, and Minhyong Kim.
\newblock Appendix and erratum to ``{M}assey products for elliptic curves of
  rank 1'' [mr2629986].
\newblock {\em J. Amer. Math. Soc.}, 24(1):281--291, 2011.

\bibitem[BR]{BakerRibet}
Matthew~H. Baker and Kenneth~A. Ribet.
\newblock Galois theory and torsion points on curves.
\newblock {\em J. Th\'eor. Nombres Bordeaux}, 15(1):11--32, 2003.
\newblock Les XXII{\`e}mes Journ{\'e}es Arithmetiques (Lille, 2001).

\bibitem[CK]{CoatesKim}
John Coates and Minhyong Kim.
\newblock Selmer varieties for curves with {CM} {J}acobians.
\newblock {\em Kyoto J. Math.}, 50(4):827--852, 2010.

\bibitem[Col]{coleman}
Robert~F. Coleman.
\newblock Dilogarithms, regulators and {$p$}-adic {$L$}-functions.
\newblock {\em Invent. Math.}, 69(2):171--208, 1982.

\bibitem[CTT]{Tzermias}
Robert~F. Coleman, Akio Tamagawa, and Pavlos Tzermias.
\newblock The cuspidal torsion packet on the {F}ermat curve.
\newblock {\em J. Reine Angew. Math.}, 496:73--81, 1998.

\bibitem[DC]{mtmueII}
Ishai Dan-Cohen.
\newblock Mixed tate motives and the unit equation {II}.
\newblock Preprint. arXiv:1510.01362.

\bibitem[DCW1]{sage1}
Ishai Dan-Cohen and Stefan Wewers.
\newblock Sage code.
\newblock \url{http://www.uni-ulm.de/mawi/rmath/mitarbeiter/wewers.html}.

\bibitem[DCW2]{CKtwo}
Ishai Dan-Cohen and Stefan Wewers.
\newblock Explicit {C}habauty-{K}im theory for the thrice punctured line in
  depth 2.
\newblock {\em Proc. Lond. Math. Soc. (3)}, 110(1):133--171, 2015.

\bibitem[DCW3]{mtmue}
Ishai Dan-Cohen and Stefan Wewers.
\newblock Mixed {T}ate motives and the unit equation.
\newblock {\em Int. Math. Res. Not. IMRN}, (17):5291--5354, 2016.

\bibitem[Del]{Deligne89}
Pierre Deligne.
\newblock Le groupe fondamental de la droite projective moins trois points.
\newblock In {\em Galois groups over {${\bf Q}$} ({B}erkeley, {CA}, 1987)},
  volume~16 of {\em Math. Sci. Res. Inst. Publ.}, pages 79--297. Springer, New
  York, 1989.

\bibitem[Fad]{Fadeev}
D.~K. Faddeev.
\newblock The group of divisor classes on some algebraic curves.
\newblock {\em Soviet Math. Dokl.}, 2:67--69, 1961.

\bibitem[Fur1]{FurushoI}
Hidekazu Furusho.
\newblock {$p$}-adic multiple zeta values. {I}. {$p$}-adic multiple
  polylogarithms and the {$p$}-adic {KZ} equation.
\newblock {\em Invent. Math.}, 155(2):253--286, 2004.

\bibitem[Fur2]{FurushoII}
Hidekazu Furusho.
\newblock {$p$}-adic multiple zeta values. {II}. {T}annakian interpretations.
\newblock {\em Amer. J. Math.}, 129(4):1105--1144, 2007.

\bibitem[Hai]{HainHigher}
Richard~M. Hain.
\newblock Higher {A}lbanese manifolds.
\newblock In {\em Hodge theory ({S}ant {C}ugat, 1985)}, volume 1246 of {\em
  Lecture Notes in Math.}, pages 84--91. Springer, Berlin, 1987.

\bibitem[Kat]{Kato}
Kazuya Kato.
\newblock Lectures on the approach to {I}wasawa theory for {H}asse-{W}eil
  {$L$}-functions via {$B_{\rm dR}$}. {I}.
\newblock In {\em Arithmetic algebraic geometry ({T}rento, 1991)}, volume 1553
  of {\em Lecture Notes in Math.}, pages 50--163. Springer, Berlin, 1993.

\bibitem[Kim1]{KimRecip}
Minhyong Kim.
\newblock Diophantine geometry and non-abelian reciprocity laws {I}.
\newblock arXiv:1312.7019.

\bibitem[Kim2]{kimi}
Minhyong Kim.
\newblock The motivic fundamental group of {$\mathbb P^1\setminus
  \{0,1,\infty\}$} and the theorem of {S}iegel.
\newblock {\em Invent. Math.}, 161(3):629--656, 2005.

\bibitem[Kim3]{kimii}
Minhyong Kim.
\newblock The unipotent {A}lbanese map and {S}elmer varieties for curves.
\newblock {\em Publ. Res. Inst. Math. Sci.}, 45(1):89--133, 2009.

\bibitem[Kim4]{KimMassey}
Minhyong Kim.
\newblock Massey products for elliptic curves of rank 1.
\newblock {\em J. Amer. Math. Soc.}, 23(3):725--747, 2010.

\bibitem[Kim5]{KimpL}
Minhyong Kim.
\newblock {$p$}-adic {$L$}-functions and {S}elmer varieties associated to
  elliptic curves with complex multiplication.
\newblock {\em Ann. of Math. (2)}, 172(1):751--759, 2010.

\bibitem[Kim6]{KimRemarks}
Minhyong Kim.
\newblock Remark on fundamental groups and effective {D}iophantine methods for
  hyperbolic curves.
\newblock In {\em Number theory, analysis and geometry}, pages 355--368.
  Springer, New York, 2012.

\bibitem[Kol]{Kolyvagin}
Victor~Alecsandrovich Kolyvagin.
\newblock On the {M}ordell-{W}eil group and the {S}hafarevich-{T}ate group of
  modular elliptic curves.
\newblock In {\em Proceedings of the {I}nternational {C}ongress of
  {M}athematicians, {V}ol.\ {I}, {II} ({K}yoto, 1990)}, pages 429--436. Math.
  Soc. Japan, Tokyo, 1991.

\bibitem[KT]{KimTamagawa}
Minhyong Kim and Akio Tamagawa.
\newblock The {$l$}-component of the unipotent {A}lbanese map.
\newblock {\em Math. Ann.}, 340(1):223--235, 2008.

\bibitem[MST]{MazurSteinTate}
Barry Mazur, William Stein, and John Tate.
\newblock Computation of {$p$}-adic heights and log convergence.
\newblock {\em Doc. Math.}, (Extra Vol.):577--614, 2006.

\bibitem[Mum]{MumfordAbVars}
David Mumford.
\newblock {\em Abelian varieties}, volume~5 of {\em Tata Institute of
  Fundamental Research Studies in Mathematics}.
\newblock Published for the Tata Institute of Fundamental Research, Bombay; by
  Hindustan Book Agency, New Delhi, 2008.
\newblock With appendices by C. P. Ramanujam and Yuri Manin, Corrected reprint
  of the second (1974) edition.

\bibitem[MvdGE]{Moonen}
Ben Moonen, Gerard van~der Geer, and Bas Edixhoven.
\newblock Abelian varieties.
\newblock \url{https://www.math.ru.nl/~bmoonen/research.html#bookabvar}.

\bibitem[NSW]{NeukirchCoh}
J{\"u}rgen Neukirch, Alexander Schmidt, and Kay Wingberg.
\newblock {\em Cohomology of number fields}, volume 323 of {\em Grundlehren der
  Mathematischen Wissenschaften [Fundamental Principles of Mathematical
  Sciences]}.
\newblock Springer-Verlag, Berlin, second edition, 2008.

\bibitem[Ols]{OlssonTowards}
Martin~C. Olsson.
\newblock Towards non-abelian {$p$}-adic {H}odge theory in the good reduction
  case.
\newblock {\em Mem. Amer. Math. Soc.}, 210(990):vi+157, 2011.

\bibitem[Ser]{Serre}
J.-P. Serre.
\newblock Local class field theory.
\newblock In {\em Algebraic {N}umber {T}heory ({P}roc. {I}nstructional {C}onf.,
  {B}righton, 1965)}, pages 128--161. Thompson, Washington, D.C., 1967.

\bibitem[SGA]{SGAIV}
{\em Groupes de monodromie en g\'eom\'etrie alg\'ebrique. {I}}.
\newblock Lecture Notes in Mathematics, Vol. 288. Springer-Verlag, Berlin-New
  York, 1972.
\newblock S\'eminaire de G\'eom\'etrie Alg\'ebrique du Bois-Marie 1967--1969
  (SGA 7 I), Dirig\'e par A. Grothendieck. Avec la collaboration de M. Raynaud
  et D. S. Rim.

\bibitem[Sil1]{SilvermanComputing}
Joseph~H. Silverman.
\newblock Computing heights on elliptic curves.
\newblock {\em Math. Comp.}, 51(183):339--358, 1988.

\bibitem[Sil2]{SilAd}
Joseph~H. Silverman.
\newblock {\em Advanced topics in the arithmetic of elliptic curves}, volume
  151 of {\em Graduate Texts in Mathematics}.
\newblock Springer-Verlag, New York, 1994.

\bibitem[Sil3]{SilArEl}
Joseph~H. Silverman.
\newblock {\em The arithmetic of elliptic curves}, volume 106 of {\em Graduate
  Texts in Mathematics}.
\newblock Springer, Dordrecht, second edition, 2009.

\bibitem[Sou]{soule}
C.~Soul{\'e}.
\newblock {$K$}-th\'eorie des anneaux d'entiers de corps de nombres et
  cohomologie \'etale.
\newblock {\em Invent. Math.}, 55(3):251--295, 1979.

\bibitem[Tat]{Tate}
J.~T. Tate.
\newblock Global class field theory.
\newblock In {\em Algebraic {N}umber {T}heory ({P}roc. {I}nstructional {C}onf.,
  {B}righton, 1965)}, pages 162--203. Thompson, Washington, D.C., 1967.

\end{thebibliography}

\bibliographystyle{alphanum}

\vfill

\noindent
{\footnotesize  J.B.: Boston University, Department of Mathematics and Statistics, 
111 Cummington Mall, Boston MA 02215, USA}

\smallskip
\noindent
{\footnotesize M.K: Mathematical Institute, University of Oxford, Woodstock Road, Oxford, OX2 6GG, and Department of Mathematics, Ewha Women's University, 52 Ewha-yeo-dae-gil, Seoul, Korea 120-750}

\smallskip
\noindent
{\footnotesize I.D.: 
Department of Mathematics, Ben-Gurion University of the Negev, P.O. Box 653, Beer-Sheva, Israel}

\smallskip
\noindent
{\footnotesize  S.W.: Institut f\"ur Reine Mathematik,
Universit\"at Ulm,
Helmholtzstrasse 18,
89081 Ulm, Germany }

\end{document}